\documentclass[12pt]{article}
\usepackage{fullpage}
\usepackage{amssymb}
\def\ad{\mathrm{ad}}

\def\b{\flat}
\def\C{\mathbb{C}}
\def\ch{\mathrm{ch}}
\def\cm{\raise.1em\hbox{$\scriptstyle\bullet$}}
\def\col{\mathrm{collar}}
\def\det{\mathrm{det}}
\def\nn{\Psi}
\def\nz{\xi}
\def\zz{\psi}
\def\ev{\mathrm{ev}}
\def\End{\mathrm{End}\,}
\def\ff{\mathrm I\!\mathrm I}
\def\FS{\mathrm{FS}}

\def\Hom{\mathrm{Hom}\,}
\def\id{\mathrm{id}}
\def\imt{\mathrm{inter}}
\def\ins{\,\lrcorner\,\,}
\def\k{\mathfrak{k}}
\def\K{\mathfrak{C}}
\def\L{\Lambda}
\def\odd{\mathrm{odd}}
\def\out{\mathrm{out}}
\def\p{\mathfrak{p}}
\def\pr{\mathrm{pr}}
\def\R{\mathbb{R}}
\def\Ric{\mathrm{Ric}}
\def\su{\mathfrak{su}}
\def\S{\mathrm{Sym}}
\def\SU{\mathbf{SU}}

\def\tr{\mathrm{tr}}
\def\U{\mathbf{U}}
\def\vol{\mathrm{vol}}
\def\Z{\mathbb{Z}}
\def\#{\sharp}
\def\<#1,#2>{\langle\,#1,\,#2\,\rangle}
\def\proof{\noindent\textbf{Proof:}\quad}
\def\pfill{\par\vskip3mm plus1mm minus1mm\noindent}
\def\qed{\ensuremath{\hfill\Box}}
\newtheorem{Lemma}{Lemma}[section]

\newtheorem{Theorem}[Lemma]{Theorem}
\newtheorem{Corollary}[Lemma]{Corollary}
\newtheorem{Definition}[Lemma]{Definition}
\newtheorem{Remark}[Lemma]{Remark}
\begin{document}
\title{About the Eta--Invariants of Berger Spheres}
\author{Gregor Weingart\footnote{
 address:\ Instituto de Matem\'aticas, Universidad Nacional Aut\'onoma
 de M\'exico, Avenida Universidad s/n, Colonia Lomas de Chamilpa, 62210
 Cuernavaca, Morelos, MEXIQUE; email: \texttt{gw@matcuer.unam.mx}.}}
\maketitle
\begin{center}
 \textbf{Abstract}
 \\[11pt]
 \parbox{412pt}{\quad\ \ 
  The integral of the top dimensional term of the multiplicative sequence
  of Pontryagin forms associated to an even formal power series is calculated
  for special Riemannian metrics on the unit ball of a hermitean vector space.
  Using this result we calculate the generating function of the reduced Dirac
  and signature $\eta$--invariants for the family of Berger metrics on the
  odd dimensional spheres.}
 \\[11pt]
\end{center}
\section{Introduction}
 Originally $\eta$--invariants of Dirac operators were introduced by
 Atiyah, Patodi and Singer in the course of establishing an index theorem
 for compact manifolds with boundary. In essence the $\eta$--invariant
 is a kind of error term arising from the geometry of the boundary,
 heuristically it can be interpreted as an expectation value for how
 many sections in the kernel and cokernel of the Dirac operator can be
 extended as $L^2$--integrable sections to the two ends of an infinite
 cylinder based on the boundary. Although the $\eta$--invariants are a
 spectral invariant of the underlying Dirac operator and vary with the
 Riemannian geometry of boundary, it appears that they somehow capture
 rather subtle $\mathbb{Q}/\Z$--valued information of the diffeomorphism
 type of the boundary manifold as exemplified by Kreck--Stolz invariants
 for $7$--manifolds. Despite their genesis as an error term $\eta$--invariants
 have thus become an object of interest of their own, whether or not the
 given manifold is a boundary.

 In this article we will focus on the calculation of the $\eta$--invariants
 of the untwisted Dirac operator and the signature operator for the family
 of Berger metrics on the odd dimensional spheres $S^{2n-1}$. Instead of
 persuing the original definition of the $\eta$--invariant as a spectral
 invariant we will solve the index formula of Atiyah--Patodi--Singer for
 the $\eta$--invariant. Needless to say this approach leads the original
 raison d'\^etre of $\eta$--invariants ad adsurdum, moreover it only allows
 us to recover the $\eta$--invariant modulo $\Z$. Interestingly the spectrum
 of the untwisted Dirac operator is known explicitly for the Berger metrics
 due to work of Hitchin \cite{ht1} in dimension $3$ and B\"ar \cite{b} in
 general. Based on these explicit spectrum calculations the $\eta$--invariant
 of the untwisted Dirac operator on the Berger spheres has been calculated
 in dimensions up to $n\,=\,16$ by Habel \cite{hab} using analytic continuation
 technics. From a much more general point of view Bechtluft--Sachs \cite{bes}
 has calculated $\eta$--invariants for low dimensional manifolds with free
 circle actions.

 \pfill
 Consider a hermitean vector space $\p$ of dimension $n$. The real part
 of the hermitean form is a scalar product $g$ on $\p$ making $\p$ a
 euclidean vector space endowed with an orthogonal complex structure
 $I:\,\p\longrightarrow\p,\,X\longmapsto iX,$ satisfying $g(IX,IY)=g(X,Y)$
 for $X,\,Y\,\in\,\p$. Rescaling the round metric on the sphere $S^{2n-1}$
 of unit vectors in $\p$ differently on the distinguished subspaces $\R IP$
 and $\{\,P,\,IP\,\}^\perp$ of $T_PS^{2n-1}\,=\,\{P\}^\perp$ defines the
 two--parameter family of Berger metrics on $S^{2n-1}$. Proportional
 Riemannian metrics share the same Levi--Civita connection and are thus
 virtually indistinguishable for the purpose of this article, so we are
 left with one geometrically relevant parameter, the ratio $T$ between
 the radii of Hopf and ``normal'' great circles. Instead of $T$ we prefer
 to use the equivalent $\rho\,:=\,T^2-1$ so that a Berger metric of
 parameter $\rho>-1$ is a Riemannian metric on $S^{2n-1}$ proportional to
 $$
  g^\rho\;\;:=\;\;g\;+\;\rho\,(\,\gamma\,\otimes\,\gamma\,)
 $$
 where $\gamma_P(X)\,:=\,g(IP,X)$ is the contact form of the natural
 CR--structure induced on the real hypersurface $S^{2n-1}$ in the flat
 K\"ahler manifold $\p$. In order to solve the index formula of
 Atiyah--Patodi--Singer for manifolds with boundary for the $\eta$--invariants
 we need to know the complementary integral of the index density over the unit
 ball $B^{2n}\,\subset\,\p$ with boundary $S^{2n-1}$. At least for the
 untwisted Dirac operator and the signature operator these index densities
 are multiplicative sequences of Pontryagin forms with complementary
 integrals given by:

 \begin{Theorem}[Special Values of Multiplicative Sequences]
 \hfill\label{mult}\break
  For every smooth metric $g^{\col,\rho}$ on the closed unit ball
  $B^{2n}\,\subset\,\p$, which is a product of the standard metric
  on $]-\varepsilon,0]$ with a Berger metric of parameter $\rho>-1$
  in a collar neighborhood $]-\varepsilon,0]\,\times\,S^{2n-1}$ of
  the boundary, the multiplicative sequence $F(\,TB^{2n},\,\nabla^{\col,
  \rho}\,)$ of Pontryagin forms associated to an even formal power series
  $F(z)=1+O(z^2)$ integrates to
  $$
   \int_{B^{2n}}F(\;TB^{2n},\,\nabla^{\col,\rho}\;)
   \;\;=\;\;
   \rho^n\,\mathrm{res}_{z=0}\left[\;\frac{F(z)^n}{z^{n+1}}\,dz\;\right]
   \;\;=\;\;
   \rho^n\,\mathrm{res}_{z=0}\left[\;\frac{(\log\,\phi)'(z)}{z^n}\,
   dz\;\right]
  $$
  where $\phi(z)=z+O(z^3)$ is the composition inverse of the formal
  power series $\frac z{F(z)}=z+O(z^3)$.
 \end{Theorem}

 \pfill
 In essence Theorem \ref{mult} is a consequence of a striking integrability
 condition enjoyed by the Berger spheres $S^{2n-1}$ thought of as the
 boundary of geodesic distance balls in $\C P^n$. Hypersurfaces in Riemannian
 manifolds satisfying this integrability condition are called permeable in
 the sequel, they are discussed in more detail in Section \ref{gh}, in
 particular we will show that all totally geodesic hypersurfaces and all
 hypersurfaces in a space form are permeable. Section \ref{trag} is devoted
 to a study of multiplicative sequences of Pontryagin forms and their
 associated (logarithmic) transgression forms, moreover we recall Hirzebruch's
 calculation of the values of multiplicative sequences of Pontryagin forms on
 $\C P^n$. The subsequent Section \ref{htra} is the technical core of this
 article, its main result Corollary \ref{logt} provides us with a closed
 formula for the logarithmic transgression form of a permeable hypersurface.
 The preceeding Lemma \ref{etf} may be of independent interest, because it
 gives a managable formula for the logarithmic transgression form of arbitrary
 hypersurfaces. The calculations for the Berger spheres proving Theorem
 \ref{mult} are relegated to the final Section \ref{bmsq}.
 
 \pfill
 Our main motivation for studying multiplicative sequences on Berger spheres
 is the relationship between $\eta$--invariants and indices of twisted Dirac
 operators as formulated in the Atiyah--Patodi--Singer Index Theorem \cite{aps}
 for manifolds with boundary. More precisely let us consider a $\Z_2$--graded
 Clifford module bundle $EM^\out$ on an even--dimensional oriented Riemannian
 manifold $M^\out$ and a totally geodesic hypersurface $M$ in $M^\out$, which
 can be realized as the boundary hypersurface $M\,:=\,\partial W$ of a compact
 subset $W\,\subset\,M^\out$ with non--trivial interior. The $\Z_2$--graded
 Clifford module bundle $EM^\out$ on $M^\out$ induces a Clifford module bundle
 $E_+M$ on $M$ and under some additional, rather severe restriction on the
 geometry of $EM^\out$ in a collar neighborhood of $M$ in $M^\out$ the Index
 Theorem of Atiyah--Patodi--Singer relates the $\eta$--invariant of the Dirac
 operator $D$ associated to $E_+M$ to the integral of the usual index density
 of the $\Z_2$--graded Clifford module bundle $EM^\out$ over $W$:
 \begin{equation}\label{apscong}
  \eta_D
  \;\;\equiv\;\;
  2\,\int_W\widehat{A}(\,TM^\out,\nabla^\out\,)\,
  \ch\,(\,EM^\out\,:\,\$M^\out\,)\qquad\textrm{mod}\;\;\Z
 \end{equation}
 Motivated by this congruence we define the reduced $\eta$--invariant of $D$
 to be any representative $\overline\eta_D\,\in\,\R$ of the class of $\eta_D$
 in $\R/\Z$ in the sense $\overline\eta_D\,\equiv\,\eta_D$ modulo $\Z$.

 In general the index density $\widehat{A}(TM^\out,\nabla^\col)
 \,\ch(EM^\out:\$M^\out)$ is not a multiplicative sequence of Pontryagin
 forms and so we can not apply Theorem \ref{mult} to calculate the interior
 integral in (\ref{apscong}). There are two well--known exceptions however,
 the untwisted Dirac operator $D^\out_\$$ on an even--dimensional spin
 manifold $M^\out$ and the signature operator $D^\out_{\mathrm{sign}}\,=\,
 d+d^*$ on an oriented manifold $M^\out$ of dimension divisible by $4$.
 The index density of the untwisted Dirac operator say is the multiplicative
 sequence of Pontryagin forms associated to the formal power series
 $\widehat{A}(z)\,:=\,\frac{\frac z2}{\sinh\,\frac z2}$, hence we obtain
 from the Atiyah--Patodi--Singer Theorem (\ref{apscong}):

 \begin{Corollary}[Eta Invariants of the Untwisted Dirac Operator]
 \hfill\label{etainv}\break
  The generating function for the reduced $\eta$--invariants $\overline{\eta}_D
  (S^{2n-1},g^\rho)$ of the untwisted Dirac operator for the Berger metric
  $g^\rho\,:=\,g\,+\,\rho\gamma\otimes\gamma$ on the odd dimensional spheres
  $S^{2n-1}$ reads:
  $$
   1\;+\;\frac12\,\sum_{n>0}\overline{\eta}_D(\;S^{2n-1},\,g^\rho\;)\,z^n
   \;\;=\;\;
   z\frac d{dz}\;\log(\;2\,\mathrm{arsinh}\,\frac{\rho\,z}2\;)
  $$
 \end{Corollary}

 \noindent
 In order to rewrite this result into a more managable formula we observe that
 the composition inverse of a solution $f$ to a differential equation $f'\,=\,
 p(f)$ is an antiderivative of $\frac1{p(z)}$. In consequence the differential
 equation $\frac d{dz}\sinh\,z\,=\,(1+\sinh^2z)^{\frac12}$ obeyed by $\sinh\,z$
 tells us
 $$
  \frac d{dz}\,\mathrm{arsinh}\,z
  \;\;=\;\;
  (\,1\;+\;z^2\,)^{-\frac12}
  \;\;=\;\;
  \sum_{k\geq 0}(\,-\frac14\,)^k\,{2k\choose k}\,z^{2k}
 $$
 via Newton's expansion $(1+z)^s\,=\,\sum_{k\geq0}{s\choose k}z^k$. Moreover
 the differential operator $z\frac d{dz}$ commutes with rescalings like
 $z\rightsquigarrow\frac{\rho\,z}2$, in this way we obtain the more managable
 formula
 \begin{equation}\label{habfor}
  1\;+\;\frac12\,\sum_{n>0}\overline{\eta}_D(\;S^{2n-1},\,g^\rho\;)\,z^n
  \;\;=\;\;
  \frac{\sum_{k\geq 0}\quad(\,-\frac1{16}\,)^k\,{2k\choose k}\,(\rho\,z)^{2k}}
  {\sum_{k\geq 0}\;\frac1{2k+1}\,(\,-\frac1{16}\,)^k\,{2k\choose k}\,
  (\rho\,z)^{2k}}
 \end{equation}
 for the $\eta$--invariants of the untwisted Dirac operator on Berger
 spheres. Needless to say the resulting power series has initial coefficients
 in accordance with the calculations done by Habel \cite{hab} in dimensions
 $n\,=\,3,\ldots,16$ mentioned before. Moreover computer experiments suggest
 that the general formula conjectured by Habel at the end of his calculations
 \begin{equation}\label{hc}
  \overline{\eta}_D(\;S^{2n-1},\,g^\rho\;)
  \;\;=\;\;
  -\,2\,(-\rho)^n\,\mathrm{res}_{x\,=\,0}
  \left[\;\frac{dx}{x^n}\,\Big(\,\sum_{l\,=\,1}^n\frac{B_l(\frac n2)}l
  \,x^{n-l}\,\Big)\,{x\,+\,\frac n2\,-\,1\choose n-1}\;\right]
 \end{equation}
 in which $B_l$ denotes the $l$--th Bernoulli polynomial, agrees with the
 formula (\ref{habfor}) at least for say $n\,\leq\,500$. The combinatorial
 complexity of the conjectural formula does not really invite a try on a
 direct proof of the equivalence of both formulas.

 Turning from the untwisted Dirac operator to the signature operator we
 consider the $\Z_2$--graded Clifford bundle $EM^\out\,=\,\L\,T^*M^\out$
 of differential forms on $M^\out$ with grading operator given by Clifford
 multiplication with the complex volume form. The induced Clifford module
 bundle $E_+M$ on $M$ can be identified with the differential forms bundle
 $\L\,T^*M$ in this case in such a way that the associated Dirac operator
 becomes $D\,=\,\Gamma(d+d^*)$, where $\Gamma$ is Clifford multiplication
 with the complex volume element of $M$. By construction $D$ commutes with
 $\Gamma$, which is essentially the Hodge $*$--isomorphism, and {\em preserves}
 the parity of differential forms, in turn $D$ becomes the direct sum of two
 operators conjugated under $\Gamma$:
 \begin{equation}\label{2op}
  D\;\;=\;\;
  [\;\Gamma\,(\,d\,+\,d^*\,)\;]^\ev
  \;\oplus\;
  [\;\Gamma\,(\,d\,+\,d^*\,)\;]^\odd
  \;\;\cong\;\;
  2\,\Gamma\,(\,d\,+\,d^*\,)^\ev
 \end{equation}
 The index density for the signature operator on the other hand is the
 multiplicative sequence of Pontryagin forms associated to the formal
 power series $L(z)\,:=\,\frac z{\tanh\,z}$ and so we conclude:

 \begin{Corollary}[Eta Invariants of the Signature Operator]
 \hfill\label{siginv}\break
  The reduced $\eta$--invariants $\overline{\eta}_{\Gamma\,(d+d^*)^\ev}
  (S^{2n-1},g^\rho)$ of the signature operator $\Gamma\,(d+d^*)^\ev$ with
  respect to the Berger metric $g^\rho\,:=\,g\,+\,\rho\gamma\otimes\gamma$
  have the following generating function
  $$
   1\;+\;
   \sum_{n>0}\overline{\eta}_{\Gamma\,(d+d^*)^\ev}(\;S^{2n-1},\,g^\rho\;)\,z^n
   \;\;=\;\;
   z\frac d{dz}\;\log(\;\mathrm{artanh}\,(\rho z)\;)
   \;\;=\;\;
   \frac{\sum_{k\geq 0}\quad(\rho z)^{2k}}
   {\sum_{k\geq 0}\;\frac1{2k+1}\,(\rho\,z)^{2k}}
  $$
  whose power series expansion is implied by the differential equation
  $\frac d{dz}\tanh\,z\,=\,1-\tanh^2z$.
 \end{Corollary}

 \pfill
 The author would like to thank Christian B\"ar, the students and the
 staff of the University of Potsdam for the hospitality enjoyed during
 various stays in Potsdam. 
\section{Geometry of Permeable Hypersurfaces}\label{gh}
 In differential geometry the Berger metrics on odd dimensional spheres
 arise naturally as the Riemannian metrics induced on the distance spheres
 $S^{2n-1}_r\,\subset \,\C P^n$ of radius $r\,\in\,]\,0,\,\frac\pi2\,[$ in
 complex projective space. The identification of Berger spheres with distance
 spheres in $\C P^n$ endow $S^{2n-1}$ with a family of second fundamental forms
 or more precisely Weingarten maps $\ff\,\in\,\Gamma(\,\End\,TS^{2n-1}\,)$,
 which satisfy the strong integrability condition characterizing permeable
 hypersurfaces. Permeable hypersurfaces as introduced in this section are
 solutions $M\,\subset\,M^\out$ to a quasilinear second order differential
 equation weaker than the second order equation $\ff\,=\,0$ characterizing
 totally geodesic hypersurfaces. The permeability of the distance spheres
 in $\C P^n$ is to be seen as a rare exception to the generic situation:
 quite probable permeable hypersurfaces do not exist in sufficiently
 wrinkled manifolds $M^\out$ due to near injectivity of the symbol map
 $\S^2T^*_pM\longrightarrow\L^3T^*_pM\otimes T_pM$ discussed below.

 \pfill
 Every hypersurface $M\subset M^\out$ in a Riemannian manifold $M^\out$ is
 naturally a Riemannian manifold itself with metric $g$ induced from the
 metric $g^\out$ of $M^\out$. The restriction of the tangent bundle $TM^\out$
 to the hypersurface $M$ is a euclidian vector bundle $\left.TM^\out\right|_M$,
 which splits orthogonally into the tangent bundle $TM$ and the normal bundle
 of $M$ in $M^\out$:
 $$
  \left.TM^\out\right|_M\;\;=\;\;TM\;\oplus\;\mathrm{Norm}\,M
 $$
 The hypersurface $M$ is called coorientable if the line bundle $\mathrm{Norm}
 \,M$ is trivial or equivalently if the monodromy of the unique metric (and
 thus flat) connection on $\mathrm{Norm}\,M$ vanishes
 $$
  \pi_1(\,M\,)\;\longrightarrow\;\Z_2,
  \qquad[\,\gamma\,]\;\longmapsto\;\mathrm{or}_M(\,[\,\gamma\,]\,)\,
  \mathrm{or}_{M^\out}(\,\iota_*[\,\gamma\,]\,)
 $$
 where $\mathrm{or}_M$ and $\mathrm{or}_{M^\out}$ are the orientation
 homomorphisms of $M$ and $M^\out$ respectively. In this article we are
 eventually interested in boundary hypersurfaces $M=\partial W$ of
 compact subsets $W\subset M^\out$ in an oriented Riemannian manifold
 $M^\out$, which are automatically oriented and cooriented by the outward
 pointing normal field $N\in\Gamma(\,\mathrm{Norm}\,M\,)$, nevertheless
 some of the arguments presented below are valid in greater generality.

 By construction the euclidian vector bundle $\left.TM^\out\right|_M$
 is endowed with two metric connections, the restriction $\left.\nabla^\out
 \right|_M$ of the Levi--Civita connection of $M^\out$ and the direct sum
 connection $\nabla$ of the Levi--Civita connection of $M$ with the unique
 metric connection on $\mathrm{Norm}\,M$. The difference between these two
 metric connections is governed by the second fundamental form $\ff$ of
 $M$ in $M^\out$ with respect to a (local) normal field $N$ defined by
 $$
  \nabla^\out_XY\;\;=\;\;\nabla_XY\;+\;\ff(\,X,\,Y\,)\,N
 $$
 for two vector fields $X,\,Y$ on $M$. Alternatively we can think of the
 second fundamental form as a vector valued $1$--form on $M$, the shape
 operator or Weingarten map $\ff:\,TM\longrightarrow TM$, by means of the
 Riemannian metric or $g(\ff_X,Y):=\ff(X,Y)$. With $\left.\nabla^\out\right|_M$
 and $\nabla$ being metric their difference is a $1$--form on $M$ with
 values in the skew--symmetric endomorphisms
 \begin{equation}\label{dsc}
  \left.\nabla^\out\right|_M\;\;=\;\;\nabla\;+\;\ff\,\wedge\,N
 \end{equation}
 on $\left.TM^\out\right|_M$, which faithfully reflects the second fundamental
 form through the definition:
 $$
  (\;\ff\,\wedge\,N\;)_XZ\;\;=\;\;
  (\;\ff_X\wedge\,N\;)\,Z\;\;:=\;\;
  g^\out(\,\ff_X,\,Z\,)\,N\;-\;g^\out(\,N,\,Z\,)\,\ff_X
 $$
 The standard formula $R^{\nabla+A}=R^\nabla+d^\nabla\!A+A^2$ becomes
 the Gau\ss--Codazzi--Mainardi equation
 \begin{equation}\label{gcm}
  \left.R^\out\right|_M
  \;\;=\;\;R\;+\;(d^\nabla\ff)\,\wedge\,N\;+\;\ff^\#\,\otimes\,\ff
 \end{equation}
 where $(\ff^\#\otimes\ff)_{X,Y}Z:=g^\out(\ff_X,Z)\ff_Y-g^\out(\ff_Y,Z)\ff_X$
 is a $2$--form on $M$ with values in the skew--symmetric endomorphisms of
 $\left.TM^\out\right|_M$. The Gau\ss--Codazzi--Mainardi equation (\ref{gcm})
 describes the restriction of the curvature tensor $R^\out$ of $M^\out$ to
 $M$ as an endomorphism--valued $2$--form in terms of the curvature tensor
 $R$ of $M$ and the second fundamental form. Alternatively we can consider
 the complete restriction of the curvature tensor $R^\out$ to $M$: 

 \begin{Definition}[Intermediate Curvature of a Hypersurface]
 \hfill\label{inter}\break
  The intermediate curvature of a hypersurface $M$ in a Riemannian manifold
  $M^\out$ is the section $R^\imt$ of the bundle of algebraic curvature
  tensors on $M$ defined by the complete restriction of the curvature
  tensor $R^\out$ to $M$ using the orthogonal projection $\pr^\perp_{TM}$
  to $TM$:
  $$
   R^\imt_{X,Y}Z
   \;\;=\;\;\pr^\perp_{TM}(\;R^\out_{X,Y}Z\;)
   \;\;=\;\;R_{X,Y}Z\;+\;(\ff^\#\otimes\ff)_{X,Y}Z
  $$
  In terms of the Gau\ss--Codazzi--Mainardi equation (\ref{gcm}) this
  definition reads $R^\imt:=R+\ff^\#\otimes\ff$.
 \end{Definition}

 \noindent
 Like every section of the bundle $\mathrm{Kr}\,TM$ of algebraic curvature
 tensors on $M$ the intermediate curvature $R^\imt$ in its covariant
 incarnation $g(R^\imt_{X,Y}Z,W)$ satisfies all symmetries of the curvature
 tensor $R$ of $M$ itself. In particular we can define the intermediate
 Ricci curvature
 $$
  \Ric^\imt(\,X,\,Y\,)
  \;\;:=\;\;
  \tr_{TM}(\;Z\;\longmapsto\;R^\imt_{Z,X}Y\;)
  \;\;=\;\;
  \sum_\mu g(\;R^\imt_{E_\mu,X}Y,\,E_\mu\;)
 $$
 as a symmetric bilinear form relating via $\Ric^\imt=\Ric+\ff^2-(\tr_g\ff)
 \,\ff$ to the actual Ricci curvature $\Ric$ of $M$. The difference
 $\Delta\Ric:=\left.\Ric^\out\right|_M-\Ric^\imt$ between the intermediate
 Ricci curvature and the restriction of the Ricci curvature $\Ric^\out$
 of $M^\out$ to $M$ provides the missing piece of information in the
 Gau\ss--Codazzi--Mainardi equation (\ref{gcm}) to describe the full
 curvature tensor $R^\out$ of $M^\out$ along $M$. In fact representation
 theory tells us that the bundle $\mathrm{Kr}\,TM^\out$ of algebraic
 curvature tensors on $M^\out$ decomposes upon restriction to $M$ as:
 $$
  \begin{array}{cccccccc}
   \left.\mathrm{Kr}\;TM^\out\right|_M
   &\;\cong\;&\mathrm{Kr}\;TM &\oplus& \L^{2,1}T^*M &\oplus&\S^2T^*M \\[2pt]
   R^\out
   &\;\widehat=\;&R^\imt&\oplus&(\,d^\nabla\ff\,)^\#&\oplus&\Delta\Ric
  \end{array}
 $$
 where $\L^{2,1}T^*M$ is the kernel of the exterior multiplication
 $\L^2T^*M\otimes T^*M\longrightarrow\L^3T^*M$. Put differently the
 difference $\Delta\Ric(X,Y)=g^\out(R^\out_{N,X}Y,N)$ parametrizes
 all second order partial derivatives of the outer metric $g^\out$
 off the hypersurface $M$ up to change of coordinates.

 \begin{Definition}[Permeable Hypersurface]
 \hfill\label{perm}\break
  Consider a hypersurface $M\subset M^\out$ in a Riemannian manifold $M^\out$
  with metric $g^\out$. The curvature tensor $R$ of the induced metric $g$ on
  $M$ considered as a $2$--form on $M$ with values in $\End\,TM$ and the second
  fundamental form $\ff$ with respect to a (local) normal field $N$ considered
  as a vector--valued $1$--form combine into a vector--valued $3$--form on $M$:
  $$
   d^\nabla(\,d^\nabla\ff\,)(X,Y,Z)
   \;\;=\;\;
   (\;R\,\ff\;)(X,Y,Z)
   \;\;:=\;\;
   R_{X,Y}\,\ff_Z\;+\;R_{Y,Z}\,\ff_X\;+\;R_{Z,X}\,\ff_Y
  $$
  The hypersurface $M$ is called permeable if $R\,\ff\in\Gamma(\L^3T^*M
  \otimes TM)$ vanishes identically, i.~e.~if the covariant exterior
  derivative $d^\nabla\ff$ of the second fundamental form is covariantly
  closed.
 \end{Definition}

 \noindent
 Of course this definition makes sense even for non--coorientable
 hypersurfaces $M\subset M^\out$, because the vanishing of $R\,\ff=0$
 is independent of the change of local normal field $N\rightsquigarrow-N$.
 Similarly it is possible to replace the curvature $R$ of $M$ in
 Definition \ref{perm} by the intermediate curvature $R^\imt$, because
 the difference $R^\imt\ff-R\,\ff=(\ff^\#\otimes\ff)\,\ff=0$ vanishes
 identically, after all its definition involves the skew--symmetrization
 of the symmetric $2$--form $g(\ff,\ff)$. The resulting alternative
 definition of permeability is in a sense more natural than Definition
 \ref{perm} in that the intermediate curvature $R^\imt$ depends only on
 the first order jet $T_pM\subset T_pM^\out$ of the hypersurface $M$ in
 $p$, whereas the curvature $R$ involves the second order jet of $M$ in
 the guise of the second fundamental form $\ff$. Permeable hypersurfaces
 are thus solutions to a quasilinear, second order differential equation,
 whose symbol in $T_pM\subset T_pM^\out$ is the kernel of the linear map
 $\S^2T^*_pM\longrightarrow\L^3T^*_pM\otimes T_pM,\,\ff\longmapsto R^\imt\ff$.
 The definition of permeability adopted in this article however has the
 advantage of being directly applicable to the calculation of
 transgressions forms.

 Recall now that a hypersurface $M\subset M^\out$ is called totally geodesic
 if every geodesic starting with a vector tangent to $M$ stays in $M$ for
 all times, equivalently $M$ is totally geodesic if and only if its second
 fundamental form $\ff=0$ vanishes identically. Clearly every totally geodesic
 hypersurface is permeable as is every hypersurface in flat space. Slightly
 more general every hypersurface in a space form $M^\out$ is permeable by
 the following argument:

 \begin{Lemma}[Permeability and Ricci Curvature]
 \hfill\label{sforms}\break
  The second fundamental form $\ff:\,TM\longrightarrow TM$ of a permeable
  hypersurface $M\subset M^\out$ commutes with the Ricci and the intermediate
  Ricci endomorphism $\Ric$ and $\Ric^\imt$ of $TM$. The necessary condition
  $[\,\Ric^\imt,\,\ff\,]\,=\,0$ for permeability is already sufficient for 
  $3$--dimensional hypersurfaces $M\subset M^\out$ in $\dim\,M^\out\,=\,4$
  or for hypersurfaces in a conformally flat Riemannian manifold $M^\out$.
  In particular all hypersurfaces in a space form are permeable.
 \end{Lemma}

 \proof
 Necessity of the condition $[\,\Ric,\,\ff\,]=0$ is a straightforward
 calculation using the definition of $R\,\ff$. The trace $\sum g(R_{X,Y}
 \ff_{E_\mu},E_\mu)\,=\,-\sum\ff(E_\mu,R_{X,Y}E_\mu)$ vanishes and so
 \begin{eqnarray*}
  \sum_\mu g(\,(R\,\ff)(E_\mu,X,Y),\,E_\mu\,)
  &=&
  \sum_\mu \left(\;g(R_{E_\mu,X}\ff_Y,\,E_\mu)\;+\;g(R_{Y,E_\mu}\ff_X,\,E_\mu)
  \;\right)
  \\
  &=&
  g((\Ric\circ\ff)\,Y,X)\,-\,g((\Ric\circ\ff)\,X,Y)
 \end{eqnarray*}
 for all $X,Y\in T_pM$. With $\Ric$ and $\ff$ being symmetric endomorphisms
 the condition $R\,\ff=0$ enforces $[\,\Ric,\,\ff\,]=0$. Replacing the
 curvature tensor $R$ by the intermediate curvature $R^\imt$ in this
 argument we get $[\,\Ric^\imt,\,\ff\,]=0$, alternatively we can
 infer $[\,\Ric^\imt,\,\ff\,]\,=\,[\,\Ric,\,\ff\,]$ from the explicit
 formula $\Ric^\imt\,=\,\Ric+\ff^2-(\tr_g\ff)\ff$.
 
 Turning from necessity to sufficiency we note that common point of the
 additional assumptions is that the intermediate curvature $R^\imt$ of
 $M$ can be written as a cross product\footnote{Sometimes called the
 ``Nomizu--Kulkarni'' product, although it is actually an instance of
 the isomorphism between the two standard presentations of a Schur functor,
 in this case the Schur functor associated to \hbox{\begin{picture}(8,8)
 \put(0,0){\line(0,1){8}}\put(4,0){\line(0,1){8}}\put(8,0){\line(0,1){8}}
 \put(0,0){\line(1,0){8}}\put(0,4){\line(1,0){8}}\put(0,8){\line(1,0){8}}
 \end{picture}}.}
 \begin{eqnarray*}
  (\,g\times h\,)(\,X,\,Y,\,Z,\,W\,)
  &:=&
  +\;g(\,X,\,Z\,)\,h(\,Y,\,W\,)\;-\;g(\,X,\,W\,)\,h(\,Y,\,Z\,)
  \\
  &&+\;g(\,Y,\,W\,)\,h(\,X,\,Z\,)\;-\;g(\,Y,\,Z\,)\,h(\,X,\,W\,)
 \end{eqnarray*}
 of the Riemannian metric $g$ and a suitable symmetric bilinear form
 $h\in\Gamma(\S^2T^*M)$. In fact in dimension $3$ the intermediate
 curvature $R^\imt$ of $M$ is a cross product with $g$ like every
 other algebraic curvature tensor. In dimensions greater than $3$
 this is no longer true, but the curvature tensor of a conformally
 flat manifold $M^\out$ can still be written in the form $R^\out=g^\out
 \times h^\out$ for some $h^\out\in\Gamma(\S^2T^*M^\out)$ so that
 $R^\imt=g\times\left.h^\out\right|_M$.

 Under the stated assumptions we can thus safely assume that $R^\imt=g\times h$
 is a cross product with $g$ for suitable $h$. The symmetric bilinear form $h$
 can be recovered from the intermediate Ricci curvature $\Ric^\imt\,=\,(2-m)h-
 (\tr_gh)g$ in dimensions $m\,:=\,\dim\,M$ different from $2$, more precisely
 $[\,\Ric^\imt,\,\ff\,]\,=\,(2-m)[\,h,\,\ff\,]$ as $g$ considered as a
 symmetric endomorphism is the identity. This equation confirms the triviality
 of Lemma \ref{sforms} in dimension $m\,=\,2$, in higher dimensions however
 the condition $[\,\Ric^\imt,\,\ff\,]\,=\,0$ for hypersurfaces $M$ with
 intermediate curvature of the form $R^\imt=g\times h$ is equivalent to
 $[\,h,\,\ff\,]=0$. Summing
 \begin{eqnarray*}
  \lefteqn{g(\,R^\imt_{X,Y}\ff_Z,\,W\,)}&&\\[2pt]
  &=&
  g(X,\ff_Z)h(Y,W)-g(X,W)h(Y,\ff_Z)-g(Y,\ff_Z)h(X,W)+g(Y,W)h(X,\ff_Z)\\[2pt]
  &=&
  \ff(X,Z)h(Y,W)-\ff(Z,Y)h(X,W)
  +g(Y,W)g((h\circ\ff)Z,X)-g(X,W)g((\ff\circ h)Y,Z)  
 \end{eqnarray*}
 cyclically over $X,\,Y,\,Z\in T_pM$ we find eventually the formula
 $$
  g(\,(R^\imt\ff)_{X,Y,Z},W\,)
  \;\;=\;\;
  g(\,[\,h,\,\ff\,]\,X,\,Y\,)\,g(\,Z,\,W\,)\;+\;
  \textrm{cyclic permutations of $X,Y,Z$}
 $$
 for all $X,\,Y,\,Z,\,W\in T_pM$, which allows us to conclude $R^\imt\ff\,=\,0$
 in case $[\,h,\,\ff\,]\,=\,0$.
 \qed
\section{Multiplicative Sequences of Pontryagin Classes}\label{trag}
 In essence multiplicative sequences are a method to construct meaningful
 ``characteristic'' classes in the de Rham cohomology of a differentiable
 manifold $M$ starting with a connection $\nabla$ on a vector bundle $VM$
 over the manifold $M$ in question. In turn the top dimensional term of
 a multiplicative sequence can be integrated over a closed manifold $M$
 yielding a differentiable invariant of $M$. Although multiplicative
 sequences may not appear particularly interesting from the more general
 Chern--Weil point of view, in this section we will focus on their
 construction and evaluation on closed manifolds without reference to
 the Chern--Weil homomorphism. The advantage of staying with multiplicative
 sequences is that their associated transgression forms can be calculated
 directly from the simpler logarithmic transgression forms, which we will 
 discuss in more detail in the next section.

 \pfill
 Recall the definition of the total Chern and Pontryagin differential forms
 associated to a connection $\nabla$ on a complex or real vector bundle
 $VM$ respectively over a manifold $M$:
 \begin{eqnarray*}
  c(\;VM,\,\nabla\;)
  &:=&
  \;\det\,\Big(\;\;\;\id\;-\;\frac{R^\nabla}{2\pi i}\;\;\;\Big)\;
  \;\;:=\;\;
  \exp\left(\;-\,\sum_{k>0}\,\frac1k\;
  \tr\,\Big(\,\frac{R^\nabla}{2\pi i}\,\Big)^k\;\right)\\
  p(\;VM,\,\nabla\;)
  &:=&
  \det^{\frac12}\Big(\;\id\,-\,\Big(\frac{R^\nabla}{2\pi}\Big)^2\;\Big)
  \;\;:=\;\;
  \exp\left(\;-\sum_{k>0}\frac1{2k}\,
  \tr\,\Big(\,\frac{R^\nabla}{2\pi}\,\Big)^{2k}\,\right)
 \end{eqnarray*}
 Evidently the homogeneous components of the total Chern form, the Chern
 forms $c_k(VM,\nabla)$, are complex valued differential forms of degree
 $2k$ on $M$, but actually they are real differential forms for every
 hermitean connection $\nabla$ on $VM$, because a symmetrized product
 of hermitean matrices is again hermitean with real trace. The second
 Bianchi identity $d^\nabla R^\nabla=0$ implies that all Chern forms
 $c_k(VM,\nabla),\,k\geq 1,$ are closed.

 The reality of the homogeneous components of the total Pontryagin form, the
 Pontryagin forms $p_k(VM,\nabla),\,k\geq 1,$ of degree $4k$, goes without
 saying and again all Pontryagin forms are closed by the second Bianchi
 identity. For the sake of a uniform treatment of multiplicative sequences
 of Chern and Pontryagin forms however we need to spoil the real definition
 of the total Pontryagin form with the introduction of some spurious
 imaginary units $i$
 $$
  p(\;VM,\,\nabla\;)
  \;\;=\;\;
  \det^{\frac12}\Big(\;\id\,+\,\Big(\frac{R^\nabla}{2\pi i}\Big)^2\;\Big)
  \;\;=\;\;
  \exp\left(\;\sum_{k>0}\frac{(-1)^{k-1}}{2k}\,
  \tr\,\Big(\,\frac{R^\nabla}{2\pi i}\,\Big)^{2k}\,\right)
 $$
 although strictly speaking this formula makes no sense unless we replace the
 real vector bundle $VM$ by its complexification $VM\otimes_\R\C$. Forgetting
 about this nuisance we get
 \begin{equation}\label{logp}
  \log\,p(\;VM,\,\nabla\;)
  \;\;=\;\;
  \sum_{k>0}\frac{(-1)^{k-1}}{2k}\,
  \tr\Big(\;\frac{R^\nabla}{2\pi i}\;\Big)^{2k}
 \end{equation}
 right from the definition of the Pontryagin form. Specifying to euclidian
 plane bundles $EM$ over $M$ we remark that the curvature of a metric
 connection $\nabla$ can be written $R^\nabla\,=\,\omega^\nabla\otimes I$
 for a (local) orthogonal complex structure $I$ on $EM$ and a suitable
 $2$--form $\omega^\nabla$ so that:
 $$
  \frac12\,\tr\Big(\;\frac{R^\nabla}{2\pi i}\;\Big)^{2k}
  \;\;=\;\;
  \frac12\,\Big(\;\frac{\omega^\nabla}{2\pi}\;\Big)^{2k}\,\tr\,(\,-iI\,)^{2k}
  \;\;=\;\;
  \Big(\;\frac{\omega^\nabla}{2\pi}\;\Big)^{2k}
 $$
 In consequence we find for the powers of the first Pontryagin form
 $p_1(EM,\nabla)\,:=\,(\frac{\omega^\nabla}{2\pi})^2$
 \begin{equation}\label{ppow}
  p_1(\;EM,\,\nabla\;)^k
  \;\;=\;\;\frac12\,\tr\Big(\;\frac{R^\nabla}{2\pi i}\;\Big)^{2k}
 \end{equation}
 and conclude that all higher Pontryagin forms $p_2(EM,\nabla),\,
 p_3(EM,\nabla),\,\ldots$ vanish due to:
 $$
  \log\,p(\;EM,\,\nabla\;)
  \;\;=\;\;
  \sum_{k>0}\frac{(-1)^{k-1}}k\,p_1(\;EM,\,\nabla\;)^k
  \;\;=\;\;
  \log(\;1\,+\,p_1(\;EM,\,\nabla\;)\;)
 $$
 The corresponding argument for Chern forms is much simpler, because
 the trace is multiplicative for endomorphisms on a complex line,
 we get an identity of differential forms
 \begin{equation}\label{cpow}
  c_1(\;LM,\,\nabla\;)^k
  \;\;=\;\;
  \left(\;\tr\,\Big(-\frac{R^\nabla}{2\pi i}\,\Big)\;\right)^k
  \;\;=\;\;
  \tr\,\Big(-\frac{R^\nabla}{2\pi i}\,\Big)^k
 \end{equation}
 for every complex line bundle $LM$. For the time being we are more interested
 in multiplicative sequences of Pontryagin forms and so we will restrict our
 discussion to the case of real vector bundles from now on, mutatis mutandis
 all arguments and definitions presented immediately translate into statements
 about multiplicative sequences of Chern forms.

 \begin{Definition}[Multiplicative Sequences of Pontryagin Forms]
 \hfill\label{mseq}\break
  The multiplicative sequence of Pontryagin forms associated to
  an even formal power series $F(z)=1+O(z^2)$ with formal logarithm
  $\log\,F(z)=\sum_{k>0}f_kz^{2k}$ is the differential form
  $$
   F(\;VM,\,\nabla\;)
   \;\;=\;\;
   \det^{\frac12}F\Big(\;\frac{R^\nabla}{2\pi i}\;\Big)
   \;\;:=\;\;
   \exp\left(\;\frac12\,\sum_{k>0}f_k\;
   \tr\Big(\frac{R^\nabla}{2\pi i}\Big)^{2k}\;\right)
  $$
  associated to a connection $\nabla$ on a real vector bundle $VM$
  over a smooth manifold $M$.
 \end{Definition}

 \noindent
 The construction of multiplicative sequences is trivially multiplicative
 under the product of (admissable) even formal power series $(F\,\tilde F)
 (VM,\nabla)=F(VM,\nabla)\,\tilde F(VM,\nabla)$. Historically however the
 name refers to the multiplicativity under the direct sum of vector bundles
 \begin{equation}\label{ism}
  F(\;VM\,\oplus\,\tilde VM,\,\nabla\,\oplus\,\tilde\nabla\;)
  \;\;=\;\;
  F(\;VM,\,\nabla\;)\,F(\;\tilde VM,\,\tilde\nabla\;)
 \end{equation}
 which follows easily from the additivity of the trace of the powers of the
 curvature
 $$
  \tr\Big(\;\frac{R^{\nabla\oplus\tilde\nabla}}{2\pi i}\;\Big)^{2k}
  \;\;=\;\;\tr\Big(\;\frac{R^\nabla}{2\pi i}\;\Big)^{2k}\;+\;
  \tr\Big(\;\frac{R^{\tilde\nabla}}{2\pi i}\;\Big)^{2k}
 $$
 of the direct sum connection $\nabla\oplus\tilde\nabla$. Equation
 (\ref{ppow}) immediately implies the characteristic property of the
 multiplicative sequence of Pontryagin forms parametrized by $F$, namely
 $$
  F(\;EM,\,\nabla\;)
  \;\;=\;\;
  \exp\Big(\;\sum_{k>0}f_k\;p_1(\;EM,\,\nabla\;)^k\;\Big)
  \;\;=\;\;
  F(\;\sqrt{p_1(\;EM,\,\nabla\;)}\;)
 $$
 for every metric connection on a euclidian plane bundle $EM$.
 Last but not least we observe
 \begin{equation}\label{ms}
  F(\;VM,\,\nabla\;)
  \;\;=\;\;
  \exp\left(\;\L^F\,\log\,p(\;VM,\,\nabla\;)\;\right)
 \end{equation}
 due to equation (\ref{logp}), where $\L^F:\,\L^{4\bullet}T^*M\longrightarrow
 \L^{4\bullet}T^*M$ multiplies homogeneous forms of degree $4k$ by $(-1)^{k-1}
 \,k\,f_k$. The right hand side of (\ref{ms}) expands into a formal power
 series in the Pontryagin forms $p_k(VM,\nabla)$ independent of the vector
 bundle $VM$. Sorting the summands of this formal power series according to
 homogeneity we get the ``multiplicative sequence'' of universal polynomials
 in the Pontryagin forms referred to in the name.

 \pfill
 Among the more interesting uses of multiplicative sequences of Pontryagin
 forms is the construction of the quantized Pontryagin forms of a real vector
 bundle $VM$, which are closely related to the Chern characters of exterior
 powers of $VM$. More precisely the quantized Pontryagin forms $P_k(VM,\nabla),
 k\geq1,$ of are the coefficients of the multiplicative sequence
 $$
  P(\;t,\,VM,\,\nabla\;)
  \;\;=\;\;
  1\;+\;t\,P_1(\;VM,\,\nabla\;)\;+\;t^2\,P_2(\;VM,\,\nabla\;)\;+\;\ldots
 $$
 of Pontryagin forms associated to the formal power series $P(t,z)\,:=\,1+
 2\,t(\,\cosh\,z-1\,)$. It is not too difficult to argue by the splitting
 principle or by the direct calculation sketched below that the quantized
 Pontryagin forms $P_k(VM,\nabla)$ for a metric connection $\nabla$ on $VM$
 vanish for $k>\frac12\dim VM$ similar to the original Pontryagin forms.
 Another similarity is
 $$
  P_r(\;VM\,\oplus\,\tilde VM,\,\nabla\,\oplus\,\tilde\nabla\;)
  \;\;=\;\;
  \sum_{s=0}^rP_s(\;VM,\,\nabla\;)\,P_{r-s}(\;\tilde VM,\,\tilde\nabla\;)
 $$
 which follows directly from the multiplicativity (\ref{ism}) of the
 sequence $P(t,VM,\nabla)$. Working out the expansion (\ref{ms}) for
 $P(t,VM,\nabla)$ explicitly we find moreover a universal formula
 $$
  P_k(\,VM,\nabla\,)
  \;\;=\;\;
  p_k(\,VM,\nabla\,)\,+\,\frac1{12}\Big(\,p_1(\,VM,\nabla\,)\,
  p_k(\,VM,\nabla\,)\,-\,(k+1)\,p_{k+1}(\,VM,\nabla\,)\,\Big)\;+\;\ldots
 $$
 for the quantized Pontryagin forms independent of the dimension of $VM$.
 The resulting congruence $P_k(VM,\nabla)\equiv p_k(VM,\nabla)$ modulo
 forms of higher degree is the main motivation for us to think of the forms
 $P_k(VM,\nabla)\,\in\,\Gamma(\L^{\geq4k}T^*M)$ as quantized Pontryagin
 forms.
 
 Motivated by a similar construction in $K$--theory \cite{a} we associate
 to every real vector bundle $VM$ with connection the formal power series
 $\sum_{d\geq 0}\tau^d\,\ch\,\L^d(VM,\nabla)$ of exterior powers of $VM$.
 This association is multiplicative in the vector bundle $VM$ in the
 sense (\ref{ism}), but does not arise directly from a multiplicative
 sequence of Pontryagin forms, because in general the Chern character
 depends on the dimension of the vector bundle, too. Ignoring this point
 for the moment we find for a metric connection $\nabla$ on a real plane
 bundle $EM$
 $$
  \sum_{d\geq0}\tau^d\ch\;\L^d(\,EM,\nabla\,)
  \;\;=\;\;
  (1+\tau e^z)\,(1+\tau e^{-z})
  \;\;=\;\;
  (1+\tau)^2\Big(1+\frac{2\tau}{(1+\tau)^2}(\cosh\,z-1)\Big)
 $$
 where $z$ is formally a root of $p_1(EM,\nabla)$. Applying the splitting
 principle we conclude
 $$
  \sum_{d\geq0}\tau^d\,\ch\;\L^d(\;VM,\,\nabla\;)
  \;\;=\;\;
  (\;1\,+\,\tau\;)^{\dim\,VM}\,P(\;\frac\tau{(1+\tau)^2},\,VM,\,\nabla\;)
 $$
 for every real vector bundle $VM$ with metric connection $\nabla$, or
 in terms of coefficients:
 $$
  \ch\;\L^d(\;VM,\,\nabla\;)
  \;\;=\;\;
  \sum_{k=0}^d{\dim\,VM-2k\choose d-k}\,P_k(\;VM,\,\nabla\;)
 $$
 In passing we note that a very similar calculation proves the formula
 $$
  \ch\;\S^d(\;VM,\,\nabla\;)
  \;\;=\;\;
  \sum_{k=0}^d{\dim\,VM+d+k-1\choose d-k}\,S_k(\;VM,\,\nabla\;)
 $$
 for the Chern character of symmetric powers of a real vector bundle
 $VM$ in terms of the coefficients $S_k(VM,\nabla),\,k\geq1,$ of the
 multiplicative sequence $S(t,VM,\nabla)$ of Pontryagin forms associated
 to the formal power series $S(t,z)=(\,1-2t(\cosh\,z-1)\,)^{-1}$. The
 coefficients $S_k(VM,\nabla)\,\in\,\Gamma(\L^{\geq4k}T^*M)$ still vanish
 in degrees less than $4k$ and are independent of the dimension of $VM$,
 but the universal formula expressing them in terms of Pontryagin forms
 is considerably more complicated than the formula for quantized Pontryagin
 forms.

 \pfill
 Eventually we want to push multiplicative sequences of Pontryagin forms,
 which result in closed differential forms $F(VM,\nabla)\in\Gamma(\L^{4
 \bullet}T^*M)$ on $M$, to multiplicative sequences of Pontryagin classes
 in de Rham cohomology. The conditio sine qua non for doing so is that
 the de Rham cohomology class $F(VM)\in H^{4\bullet}_{\mathrm{dR}}(M)$
 represented by the closed differential form $F(VM,\nabla)\in\Gamma
 (\L^{4\bullet}T^*M)$ is independent of the connection $\nabla$. En nuce
 this independence relies on the existence of transgression forms, explicit
 solutions $\mathrm{Trans}\,F(VM,\nabla^0,\nabla^1)$ to the transgression
 problem to make the difference between two representative closed forms exact:
 \begin{equation}\label{tp}
  F(\;VM,\,\nabla^1\;)\;-\;F(\;VM,\,\nabla^0\,)
  \;\;=\;\;
  d\,(\mathrm{Trans}\,F)(VM,\nabla^0,\nabla^1)
 \end{equation}
 Needless to say these transgression forms are to live a life of their own,
 besides our use of them to calculate the values of multiplicative sequences
 on compact manifolds with boundary they lead to the definition of secondary
 characteristic classes of flat vector bundles. Although not completely
 straightforward we begin our discussion of transgression forms with:

 \begin{Definition}[Logarithmic Transgression Form]
 \hfill\label{logform}\break
  The logarithmic transgression form associated to a multiplicative
  sequence of Pontryagin classes parametrized by an even formal power
  series $F(z)=1+O(z^2)$ with logarithm $\log\,F(z)=\sum_{k>0}f_kz^{2k}$
  at a tangent vector $(\nabla,\dot\nabla)$ to the space of connections
  on a (real) vector bundle $VM$ over a differentiable manifold $M$ is
  the following differential form on $M$:
  $$
   \delta\,(\log\,F)(\;VM,\,\nabla,\,\dot\nabla\;)
   \;\;:=\;\;
   \sum_{k>0}\;k\,f_k\;\tr\Big(\;\frac{\dot\nabla}{2\pi i}\,
   \Big(\frac{R^\nabla}{2\pi i}\Big)^{2k-1}\;\Big)
  $$
 \end{Definition}

 \pfill
 The logarithmic transgression form describes the logarithmic derivative
 of the multiplicative sequence $\nabla\longmapsto F(\,VM,\,\nabla\,)$ along
 a curve curve $t\mapsto\nabla^t$ in the space of connections through:
 \begin{equation}\label{chtr}
  \frac d{dt}\;\log\,F\,(\;VM,\,\nabla^t\;)
  \;\;=\;\;
  d\;\left(\delta\,(\log\,F)(\;VM,\,\nabla^t,\,
  \frac d{dt}\nabla^t\;)\right)
 \end{equation}
 Upon integration over $[0,t]$ and exponentiation in the exterior algebra this
 equation becomes:
 \begin{equation}\label{this}
  F(\;VM,\,\nabla^t\;)
  \;\;=\;\;
  F(\;VM,\,\nabla^0\;)\;\exp\left(\;d\int_0^t\delta\,(\log\,F)
  (\,VM,\nabla^\tau,\frac d{d\tau}\nabla^\tau\,)\,d\tau\;\right)
 \end{equation}
 Taking the derivative of this equation in $t$ and integrating once more
 over $[0,1]$ we see that
 \begin{eqnarray}
  \lefteqn{(\,\mathrm{trans}\,F\,)(\;VM,\,\nabla^0,\,\nabla^1\;)\;\;:=}
  \qquad&&\label{deft}\\[3pt] &&
  \int_0^1\exp\left(\;\int_0^td\,\delta\,(\log\,F)
  (\,VM,\nabla^\tau,\frac d{d\tau}\nabla^\tau\,)\,d\tau\;\right)
  \;\delta\,(\log\,F)(\,VM,\nabla^t,\frac d{dt}\nabla^t\,)\;dt\nonumber
 \end{eqnarray}
 solves the following reformulation of the original transgression
 problem (\ref{tp}):
 \begin{equation}\label{teq}
  F(\;VM,\,\nabla^1\;)\;-\;F(\;VM,\,\nabla^0\;)
  \;\;=\;\;
  F(\;VM,\,\nabla^0\;)\;d\,(\mathrm{trans}\,F)(\;VM,\,\nabla^0,\,\nabla^1\;)
 \end{equation}
 Given the rather ardeous formula (\ref{deft}) for the transgression form it
 may be hard to believe that it can be of any use besides showing that the
 de Rham cohomology class $F(VM)$ represented by $F(VM,\nabla)$ is well
 defined. The importance of formula (\ref{deft}) however lies in the fact
 that the solution to the modified transgression problem (\ref{teq}) is
 expressed solely in terms of the time dependence of the logarithmic
 transgression form $\delta(\log\,F)(VM,\nabla^t,\frac{d}{dt}\nabla^t)$.
 In the explicit calculation of transgression forms in Section \ref{htra}
 this becomes a crucial advantage.

 \begin{Remark}[Multiplicative Sequences of Chern Forms]
 \hfill\label{cseq}\break
  Similarly to multiplicative sequences of Pontryagin forms (or classes)
  we can define multiplicative sequences of Chern forms for complex vector
  bundles $VM$ over a manifold $M$. The differential form associated to a
  formal power series $F(z)=1+O(z)$ with formal logarithm $\log\,F(z)=
  \sum_{k>0}f_kz^k$ and a connection $\nabla$ on a complex vector bundle
  $VM$ over $M$ reads
  $$
   F(\;VM,\,\nabla\;)
   \;\;=\;\;
   \det\,F\Big(-\frac{R^\nabla}{2\pi i}\;\Big)
   \;\;:=\;\;
   \exp\left(\;\sum_{k>0}f_k\;
   \tr\Big(-\frac{R^\nabla}{2\pi i}\Big)^k\;\right)
  $$
  while the corresponding logarithmic transgression form is defined as:
  $$
   \delta\,(\log\,F)(\;VM,\,\nabla,\,\dot\nabla\;)
   \;\;:=\;\;
   \sum_{k>0}\;k\,f_k\;\tr\Big(\;-\frac{\dot\nabla}{2\pi i}\,
   \Big(-\frac{R^\nabla}{2\pi i}\Big)^{k-1}\;\Big)
  $$
  Quantized Chern classes $C_k(VM,\nabla),\,k\geq 1,$ are defined as well
  as the coefficients of the multiplicative sequence $C(t,VM,\nabla)$
  associated to the power series $C(t,z)\,:=\,1+t(e^z-1)$.
 \end{Remark}

 \pfill
 On compact manifolds $M$ of dimension divisible by $4$ we can integrate
 the top term of the multiplicative sequence $F(TM,\nabla)$ of Pontryagin
 forms over $M$. For closed manifolds the resulting value $\<F(TM,\nabla),
 [M]>$ of the multiplicative sequence on $M$ only depends on the de Rham
 cohomology class $F(TM)$ of the differential form $F(TM,\nabla)$
 and is thus independent of the connection $\nabla$, for compact manifolds
 with boundary things are more difficult and the dependence of $\<F(TM,\nabla),
 [M]>$ on $\nabla$ inevitably involves the values of transgression forms.
 In general it is quite difficult to calculate the values of multiplicative
 sequences or transgression forms, for complex projective spaces $\C P^n$
 however the problem can be solved completely. The complex tangent bundle
 of the complex projective space is isomorphic to
 $$
  T^{1,0}\C P^n\;\;=\;\;\Hom(\;\mathcal{O}(-1),\,\C^{n+1}/
  _{\displaystyle\mathcal{O}(-1)}\;)
  \;\;=\;\;(\,n+1\,)\,\mathcal{O}(1)\;-\;\C
 $$
 where $\C$ and $\C^{n+1}$ denote the trivial complex vector bundles over
 $\C P^n$ of rank $1$ and $n+1$ and $\mathcal{O}(-1)$ is the tautological
 complex line bundle with dual ``hyperplane'' bundle $\mathcal{O}(1)$.
 Forgetting the complex structure we get an isomorphism of real vector
 bundles
 \begin{equation}\label{dsum}
  T\,\C P^n\,\oplus\,\R^2
  \;\;=\;\;
  (\,n+1\,)\,\mathcal{O}(-1)^\R
 \end{equation}
 which results in the equality $F(T\C P^n)=F(\mathcal{O}(-1)^\R)^{n+1}$
 for every multiplicative sequence of Pontryagin classes as trivial or
 more generally flat vector bundles $VM$ only contribute a factor $F(VM)=1$.
 In order to calculate the first Pontryagin class of the real plane bundle
 $\mathcal{O}(-1)^\R$ we consider its complexification $\mathcal{O}(-1)^\R
 \otimes_\R\C\cong\mathcal{O}(-1)\oplus\mathcal{O}(1)$, which tells us
 $$
  p_1(\;\mathcal{O}(-1)^\R\;)
  \;\;=\;\;
  -\,c_2(\;\mathcal{O}(-1)\oplus\mathcal{O}(1)\;)
  \;\;=\;\;
  -\,c_1(\;\mathcal{O}(-1)\;)\,c_1(\;\mathcal{O}(1)\;)
  \;\;=\;\;\Big(\frac{\omega^\FS}\pi\Big)^2
 $$
 where $\omega^\FS$ is the K\"ahler form of the Fubini--Study metric on
 $\C P^n$, because the first Chern form of the complex line bundles
 $\mathcal{O}(k)$ with respect to a Fubini--Study type connection reads:
 $$
  c_1(\;\mathcal{O}(\,k\,),\,\nabla^\FS\;)\;\;=\;\;k\,\frac{\omega^\FS}\pi
 $$
 All in all our consideration starting with the isomorphism (\ref{dsum})
 imply the following equality
 $$
  F(\;T\C P^n,\,\nabla^\FS\;)\;\;=\;\;F(\;\frac{\omega^\FS}\pi\;)^{n+1}
 $$
 in de Rham cohomology, which must be true already on the level of differential
 forms, because both sides are left invariant and thus parallel on the
 symmetric space $\C P^n$. Using the volume integral $\int_{\C P^n}\frac1{n!}
 (\omega^\FS)^n=\frac{\pi^n}{n!}$ of the complex projective space we find
 $\<(\frac{\omega^\FS}\pi)^n,[\C P^n]>=1$ or:
 \begin{equation}\label{creq}
  \<\;F(\;T\C P^n\;),\,[\,\C P^n\,]>
  \;\;=\;\;
  \mathrm{res}_{z=0}\left[\frac{F(z)^{n+1}}{z^{n+1}}dz\right]
 \end{equation}
 In this way we have reduced the problem of calculating the values
 $\<F(T\C P^n),[\C P^n]>$ of the multiplicative sequence $F$ on complex
 projective spaces to the combinatorial problem of finding the residues
 on the right hand side of equation (\ref{creq}). Introducing the odd
 formal power series $\phi(z)=z+O(z^3)$ explicitly as the (composition)
 inverse of $\frac z{F(z)}=z+O(z^3)$ or implicitly by $F(\phi(z))\,=\,
 \frac{\phi(z)}z$ we can apply the transformation rule for residues to get:
 $$
  \mathrm{res}_{z=0}\left[\frac{F(z)^{n+1}}{z^{n+1}}\,dz\right]
  \;\;=\;\;
  \mathrm{res}_{z=0}\left[\frac{F(\phi(z))^{n+1}}{\phi(z)^{n+1}}
  \,\phi'(z)\,dz\right]
  \;\;=\;\;
  \mathrm{res}_{z=0}\left[\frac{\phi'(z)}{z^{n+1}}\,dz\right]
 $$

 \begin{Theorem}[Hirzebruch's Proof of the Signature Theorem \cite{hir}]
 \hfill\label{hproof}\break
  Consider a multiplicative sequence of Pontryagin classes associated
  to an even formal power series $F(z)=1+O(z^2)$. The odd formal power
  series $\phi(z)=z+O(z^3)$ defined as the composition inverse of the
  odd power series $\frac z{F(z)}=z+O(z^3)$ or implicity by $F(\phi(z))
  =\frac{\phi(z)}z$ is the generating power series for the values of
  the multiplicative sequence associated to $F$
  $$
   \phi'(z)\;\;=\;\;1\;+\;\sum_{k>0}\<F(\;T\C P^{2k}\;),[\;\C P^{2k}\;]>z^{2k}
  $$
  on the complex projective spaces. Note that $\<F(T\C P^k),[\C P^k]>=0$
  for odd $k$ by definition.
 \end{Theorem}

 \pfill
 Using Theorem \ref{hproof} we can construct multiplicative sequences of
 Pontryagin classes assuming arbitrarily prescribed values on the even
 and vanishing on the odd complex projective spaces. Collecting these
 values in the formal power series $\phi'(z):=1+\sum_{k>0}\phi'_kz^{2k}$
 we simply choose the multiplicative sequence parametrized by the quotient
 $F(z)=\frac z{f(z)}$, where $f$ is the unique formal power series solution
 to the differential equation $f'=\frac1{\phi'(f)}$ with initial value
 $f(0)=0$. In particular there is a unique multiplicative sequence which
 takes the values $1$ and $0$ on even and odd complex projective spaces
 respectively, namely the sequence associated to
 $$
  L(z)\;\;=\;\;\frac z{\tanh\,z}
 $$
 because $l(z)=\tanh\,z$ is the unique solution to the differential equation
 $l'=1-l^2$ with initial value $l(0)=0$. On the other hand the even and odd
 complex projective spaces have signature $1$ and $0$ respectively and generate
 the rational oriented cobordism ring due to a result of Thom. With signature
 being an oriented cobordism invariant Hirzebruch concluded
 $$
  \mathrm{sign}\;M\;\;=\;\;\<L(TM),[M]>
 $$
 for every closed manifold $M$. Instead of pursuing these ideas further and
 start a discussion of general rational oriented cobordism invariants or
 ``genera'' we want to use Theorem \ref{hproof} to calculate the values of
 a particular genus, the $\widehat{A}$--genus, parametrized by
 $$
  \widehat{A}(z)\;\;=\;\;\frac{\frac z2}{\sinh\,\frac z2}
 $$
 on the complex projective spaces. With $\widehat{a}(z)=2\sinh\frac z2$
 being a solution to the differential equation $\widehat{a}'(z)=(1+\frac14
 \widehat{a}^2(z))^{\frac12}$ we immediately find the generating formal
 power series
 \begin{equation}\label{ahat}
  1\;+\;\sum_{k>0}\<\widehat{A}(T\C P^{2k}),[\C P^{2k}]>\,z^{2k}
  \;\;=\;\;
  \Big(\;1+\frac{z^2}4\;\Big)^{-\frac12}
  \;\;=\;\;
  \sum_{k\geq 0}\Big(-\frac1{16}\Big)^k\,{2k\choose k}\,z^{2k}
 \end{equation}
 by Newton's power series expansion $(1+z)^e=\sum_{k\geq 0}{e\choose k}z^k$
 for $e=-\frac12$ and ${-\frac12\choose k}=(-\frac14)^k{2k\choose k}$.
\section{Solution of the Hypersurface Transgression Problem}\label{htra}
 On closed manifolds the evaluation of multiplicative sequences of
 Pontryagin forms factorizes over the Pontryagin classes in de Rham
 cohomology. Similar computations on manifolds with boundary involve
 additional boundary contributions in form of the integrals of associated
 transgression forms. In this section we want to study Definition
 \ref{logform} of the logarithmic transgression form associated to a
 multiplicative sequence in more detail for the standard transgression
 problem associated to a (boundary) hypersurface $M\,\subset\,M^\out$. In
 particular the formula obtained in Corollary \ref{logt} for the logarithmic
 transgression form of a permeable hypersurface $M\,\subset\,M^\out$ only
 depends on the $3$--form $g(d^\nabla\ff,\ff)\,\in\,\Gamma(\L^3T^*M)$.

 \pfill
 In order to describe the standard transgression problem associated to a
 hypersurface let us consider a compact subset $W\subset M^\out$ of an
 oriented manifold $M^\out$ with smooth boundary $M\,=\,\partial W$ and
 fix a collar neighborhood $]-\varepsilon,\varepsilon[\times M\supset M$
 of $M$ such that the coordinate $r$ of $]-\varepsilon,\varepsilon[$ is
 zero on $M$ and negative on the interior of $W$. A Riemannian metric
 $g^\col$ on $M^\out$ is called a collar metric, if it restricts on the
 collar $]-\varepsilon,\varepsilon[\,\times M$ to the product
 $$
  \left.g^\col\right|_{]-\varepsilon,\varepsilon[\,\times M}
  \;\;=\;\;
  dr\otimes dr\;+\;g
 $$
 of the standard metric on $]-\varepsilon,\varepsilon[$ with a Riemannian
 metric $g$ on $M$, the vector field $N\,:=\,\frac{\partial}{\partial r}$
 is thus the outward pointing normal vector field to $M$. In the context
 of Theorem \ref{mult} we are interested in the value of a multiplicative
 sequence of Pontryagin forms
 \begin{equation}\label{fint}
  \<F(\;TM^\out,\,\nabla^\col\;),[\,W\,]>
  \;\;:=\;\;
  \int_WF(\;TM^\out,\,\nabla^\col\;)
 \end{equation}
 for the Levi--Civita connection $\nabla^\col$ of a collar metric $g^\col$.
 An exact evaluation of the integral (\ref{fint}) may be impossible due
 to lack of precise information concerning the collar metric $g^\col$,
 however in favorable situations there might be another Riemannian metric
 $g^\out$ on $M^\out$ such that the calculation of the integral $\<F(TM^\out,
 \nabla^\out),[\,W\,]>$ for the Levi--Civita connection of $g^\out$ is
 feasible. In this case Stokes' Theorem allows us to convert the difference
 \begin{equation}\label{tokes}
  \int_WF(TM^\out,\nabla^\out)
  \;-\;
  \int_WF(TM^\out,\nabla^\col)
  \;\;=\;\;
  \int_M\mathrm{Trans}\,F
  (TM^\out,\nabla^\col,\nabla^\out)
 \end{equation}
 into a boundary term involving the transgression form $(\mathrm{Trans}\,F)
 (TM^\out,\nabla^\col,\nabla^\out)$. In turn the integration over $M$ only
 depends on the restriction of the transgression form
 $$
  \left.(\,\mathrm{Trans}\,F\,)
  (\,TM^\out,\nabla^\col,\nabla^\out\,)\right|_M
  \;\;=\;\;
  F(\left.TM^\out\right|_M,\nabla)\,
  (\,\mathrm{trans}\,F\,)(\left.TM^\out\right|_M,\nabla,
  \left.\nabla^\out\right|_M)
 $$
 to the hypersurface $M\subset M^\out$, which agrees by naturality with the
 transgression form for the restricted connections $\nabla=\left.\nabla^\col
 \right|_M$ and $\left.\nabla^\out\right|_M$ on $\left.TM^\out\right|_M$.
 Of course the restricted connection $\nabla$ is just the Levi--Civita
 connection for the metric $g$ on $M$ extended trivially to a metric
 connection on $\left.TM^\out\right|_M\cong TM\oplus\mathrm{Norm}\,M$
 so that $F(\left.TM^\out\right|_M,\nabla)\,=\,F(TM,\nabla)$. According
 to Section \ref{gh} the homotopy $t\longmapsto t\nabla^\out+(1-t)\nabla^\col$
 between the connections $\nabla^\col$ and $\nabla^\out$ on $TM^\out$ restricts
 to the linear interpolation $t\longmapsto\nabla+t\,\ff\wedge N$ between
 $\nabla$ and $\left.\nabla^\out\right|_M$ on $\left.TM^\out\right|_M$,
 in consequence the associated transgression form depends crucially on
 the second fundamental form $\ff$ of $M\subset M^\out$. Summarizing
 these considerations we define:

 \begin{Definition}[Logarithmic Transgression Form of a Hypersurface]
 \hfill\label{thyp}\break
  The logarithmic transgression form of a hypersurface $M\subset M^\out$ in
  a Riemannian mani\-fold $M^\out$ with metric $g^\out$ is the logarithmic
  transgression form of (the derivative of) the linear interpolation
  $\nabla^t\,=\,\nabla+t\,\ff\wedge N$ between the metric connections
  $\nabla$ and $\nabla+\ff\wedge N$ on $\left.TM^\out\right|_M$:
  $$
   \delta\,(\log\,F) 
   (\,\left.TM^\out\right|_M,\nabla^t,\ff\wedge N\,)
   \;\;=\;\;
   \sum_{k>0}k\,f_k\,\tr\left(\frac{\ff\wedge N}{2\pi i}\,
   \Big(\frac{R^t}{2\pi i}\Big)^{2k-1}\right)
  $$
  In this formula the curvature $R^t$ of $\nabla^t$ is given by the
  Gau\ss--Codazzi--Mainardi equation (\ref{gcm}):
  $$
   R^t
   \;\;=\;\;R\;+\;t\,d^\nabla\ff\,\wedge\,N\;+\;t^2\,\ff^\#\,\otimes\,\ff
  $$
 \end{Definition}

 \noindent
 Although this definition of the logarithmic transgression form is useful for
 theoretical considerations, it is not suitable for actual calculations,
 because its time dependence involves the time dependent curvature $R^t$ of
 the connection $\nabla^t$ on $\left.TM^\out\right|_M$. For the rest of this
 section we want to study the definition of the logarithmic transgression form
 in more detail in order to derive a more convenient expansion for $\delta(\log
 \,F)(\,\left.TM^\out\right|_M,\nabla^t,\ff\wedge N\,)$. Let us begin this
 endeavour with the definition of the basic geometric differential forms
 \begin{equation}\label{forms}
  \nz\;\,:=\;\,\sum_{k\geq0}g(\,d^\nabla\ff,\,R^k\,\ff\,)
  \qquad
  \zz\;\,:=\;\,\sum_{k\geq0}g(\,\ff,\,R^k\,\ff\,)
  \qquad
  \nn\;\,:=\;\,\sum_{k\geq0}g(\,d^\nabla\ff,\,R^k\,d^\nabla\ff\,)
 \end{equation}
 encoding the geometry of the hypersurface $M\,\subset\,M^\out$. Evidently
 $\nz\,\in\,\Gamma(\L^\odd T^*M)$ is an odd differential form with non--zero
 homogeneous components $\nz_3,\,\nz_5,\,\ldots$ living in degrees at least
 $3$. Similarly $\zz$ and $\nn$ are differential forms with non--zero
 homogeneous components $\zz_4,\,\zz_8,\,\ldots$ and $\nn_4,\,\nn_8,\,
 \ldots$ concentrated in positive degrees divisible by $4$, because say
 $$
  g(\;\ff,\,R^k\,\ff\;)
  \;\;=\;\;(-1)^k\,g(\;R^k\,\ff,\,\ff\;)
  \;\;=\;\;-\,(-1)^k\,g(\;\ff,\,R^k\,\ff\;)
 $$
 where the stray sign change in the last equality is caused by passing the
 odd vector valued form $\ff$ past the odd form $R^k\,\ff$. A similar argument
 implies that the last possible variant
 $$
  i^{N+1}\nz
  \;\;:=\;\;\sum_{k\geq 0}g(\,\ff,\,R^k\,d^\nabla\ff\,)
  \;\;=\;\;\sum_{k\geq 0}(-1)^k\,g(\,d^\nabla\ff,\,R^k\,\ff\,)
 $$
 in the definition (\ref{forms}) of the basic geometric differential forms
 results in a differential form $i^{N+1}\nz$ on $M$, which is the image of
 $\nz$ under the operator $i^{N+1}:\,\L^\odd T^*M\longrightarrow\L^\odd T^*M$,
 which multiplies a homogeneous odd form of degree $r$ by $i^{r+1}$. Using
 the second Bianchi identity $d^\nabla R\,=\,0$ and the identity $d^\nabla
 (\,d^\nabla\ff\,)\,=\,R\,\ff$ for the metric connection $\nabla$ on $TM$
 we readily find the exterior differential system satisfied by the
 differential forms $\nz$, $\zz$ and $\nn$
 \begin{eqnarray}
  d\nz
  &=&
  \sum_{k\geq 0}\Big(\;
  g(\,R\,\ff,\,R^k\,\ff\,)\;+\;g(\,d^\nabla\ff,\,R^k\,d^\nabla\ff\,)\;\Big)
  \quad\;=\;\;
  \nn\;-\;\zz
  \nonumber
  \\
  d\zz
  &=&
  \sum_{k\geq 0}
  \Big(\;g(\,d^\nabla\ff,\,R^k\,\ff\,)
  \;-\;g(\,\ff,\,R^k\,d^\nabla\ff\,)\;\Big)
  \qquad\;=\;\;
  \nz\;-\;i^{N+1}\nz
  \label{eds}
  \\
  d\nn
  &=&
  \sum_{k\geq 0}
  \Big(\;g(\,R\,\ff,\,R^k\,d^\nabla\ff\,)
  \;+\;g(\,d^\nabla\ff,\,R^k\,R\,\ff\,)\;\Big)
  \;\;=\;\;
  \nz\;-\;i^{N+1}\nz
  \nonumber
 \end{eqnarray}
 where the missing summands care for themselves. In terms of the homogeneous
 components of the differential forms $\nz$, $\zz$ and $\nn$ these equations
 can be written schematically:
 $$
  \R\,\nz_3
  \;\stackrel d\longrightarrow\;\R\,\nn_4\,\oplus\,\R\,\zz_4
  \;\stackrel d\longrightarrow\;\R\,\nz_5
  \qquad
  \R\,\nz_7
  \;\stackrel d\longrightarrow\;\R\,\nn_8\,\oplus\,\R\,\zz_8
  \;\stackrel d\longrightarrow\;\R\,\nz_9
  \qquad
  \ldots
 $$

 \pfill
 Perhaps surprisingly the main problem in understanding the logarithmic
 transgression form of a hypersurface is to express the time dependent
 variants of the geometric differential forms
 $$
  \nz^t\;\,:=\;\,\sum_{k\geq0}g(\,d^\nabla\ff,(R^t)^k\,\ff\,)
  \qquad
  \zz^t\;\,:=\;\,\sum_{k\geq0}g(\,\ff,(R^t)^k\,\ff\,)
  \qquad
  \nn^t\;\,:=\;\,\sum_{k\geq0}g(\,d^\nabla\ff,(R^t)^k\,d^\nabla\ff\,)
 $$
 in terms of the basic geometric forms $\nz$, $\zz$ and $\nn$. In passing
 we note that the differential forms $\nz^t$, $\zz^t$ and $\nn^t$ satisfy
 the same exterior differential system (\ref{eds}) as the basic geometric
 forms $\nz$, $\zz$ and $\nn$, in the argument given above we need only
 observe $d^\nabla\ff\,=\,d^{\nabla^t}\ff$ and replace the metric connection
 $\nabla$ by the metric connection $\nabla^t$. Moreover the possibly non--zero
 homogeneous components of the differential forms $\zz^t$ and $\nn^t$ are
 concentrated in positive degrees divisible by $4$ as before, while $\nz^t$
 is an odd differential form with vanishing component in degree $1$. In
 difference to the basic geometric differential forms $\nz$, $\zz$ and
 $\nn$ however the time dependence of $\nz^t$, $\zz^t$ and $\nn^t$ allows
 us study the ordinary differential equation
 \begin{eqnarray*}
  \frac d{dt}\;\nz^t
  &=&
  -\;2\,t\,\nz^t\,\wedge\,(\;\zz^t\;+\;\nn^t\;)
  \\
  \frac d{dt}\;\zz^t
  &=&
  -\;2\,t\,\zz^t\,\wedge\,\zz^t\;+\;2\,t\,\nz^t\,\wedge\,(\,i^{N+1}\nz^t\,)
  \\
  \frac d{dt}\;\nn^t
  &=&
  -\;2\,t\,\nn^t\,\wedge\,\nn^t\;-\;2\,t\,\nz^t\,\wedge\,(\,i^{N+1}\nz^t\,)
 \end{eqnarray*}
 satisfied by the triple $(\nz^t,\,\zz^t,\,\nn^t)$, which is a direct
 consequence of the corresponding equation for the time dependent curvature
 $R^t\,=\,R\,+\,t\,d^\nabla\ff\wedge N\,+\,t^2\,\ff^\#\otimes\ff$. In fact
 we find
 \begin{eqnarray*}
  \frac d{dt}\;\nz^t
  &=&
  \sum_{k,\,l\geq0}
  g\Big(\;d^\nabla\ff,\,(R^t)^k\,
  \Big(\,(d^\nabla\ff)^\#\,\wedge\,N\,-\,N^\#\,\otimes\,(d^\nabla\ff)
  \,+\,2\,t\,\ff^\#\,\otimes\,\ff\,\Big)\,(R^t)^l\,\ff\;\Big)
  \\
  &=&
  \sum_{k,\,l\geq0}\Big(\,
  g(\,d^\nabla\ff,\,(R^t)^k\,N\,)\wedge g(\,d^\nabla\ff,\,(R^t)^l\,\ff\,)
  -
  g(\,d^\nabla\ff,\,(R^t)^k\,d^\nabla\ff\,)\wedge g(\,N,\,(R^t)^l\,\ff\,)
  \,\Big)
  \\[-7pt]
  &&
  \qquad\quad-\;2\,t\,\sum_{k,\,l\geq0}
  g(\,d^\nabla\ff,\,(R^t)^k\,\ff\,)\,\wedge\,g(\,\ff,\,(R^t)^l\,\ff\,)
  \\[2pt]
  &=&
  -\;2\,t\,\nn^t\,\wedge\,\nz^t\;-\;2\,t\,\nz^t\,\wedge\,\zz^t
 \end{eqnarray*}
 using $R^t\,N\,=\,-\,t\,d^\nabla\ff$, replacing the leftmost $d^\nabla\ff$
 by $\ff$ and the rightmost $\ff$ by $d^\nabla\ff$ respectively in the second
 line provides the analoguous differential equations for $\zz^t$ and $\nn^t$.

 \begin{Lemma}[Solution of Special Ordinary Differential Equation]
 \hfill\label{etf}\break
  Consider the exterior algebra $\L\,T^*$ of alternating forms on a vector
  space $T$ and the subspace $\L^{4\bullet}_\circ T^*\,\subset\,\L\,T^*$ of
  forms concentrated in degrees divisible by $4$ without constant term. Moreover
  let $i^{N+1}:\,\L^\odd T^*\longrightarrow\L^\odd T^*$ be the operator on
  the subspace $\L^\odd T^*\,\subset\,\L\,T^*$ of odd forms on $T$ multiplying
  a homogeneous odd form of degree $r$ by $i^{r+1}$. For every triple of forms
  $(\nz,\,\zz,\,\nn)\,\in\,\L^\odd T^*\oplus\L^{4\bullet}_\circ T^*\oplus
  \L^{4\bullet}_\circ T^*$ the formal power series $(\nz^t,\,\zz^t,\,\nn^t)$
  in $t$ with
  \begin{eqnarray*}
   \nz^t
   &:=&
   \nz\,\wedge\,(\;1\;+\;t^2\,\zz\;)^{-1}\,\wedge\,(\;1\;+\;t^2\,\nn\;)^{-1}
   \\
   \zz^t
   &:=&
   \zz\,\wedge\,(\;1\;+\;t^2\,\zz\;)^{-1}
   \;+\;t^2\,\nz\,\wedge\,(\,i^{N+1}\nz\,)
   \,\wedge\,(\;1\;+\;t^2\,\zz\;)^{-2}\,\wedge\,(\;1\;+\;t^2\,\nn\;)^{-1}
   \\
   \nn^t
   &:=&
   \nn\,\wedge\,(\;1\;+\;t^2\,\nn\;)^{-1}
   \;-\;t^2\,\nz\,\wedge\,(\,i^{N+1}\nz\,)
   \,\wedge\,(\;1\;+\;t^2\,\zz\;)^{-1}\,\wedge\,(\;1\;+\;t^2\,\nn\;)^{-2}
  \end{eqnarray*}
  is the unique solution with initial value $(\nz,\,\zz,\,\nn)$ in $0$
  to the ordinary differential equation:
  \begin{eqnarray*}
   \frac d{dt}\;\nz^t
   &=&
   -\;2\,t\,\nz^t\,\wedge\,(\;\zz^t\;+\;\nn^t\;)
   \\
   \frac d{dt}\;\zz^t
   &=&
   -\;2\,t\,\zz^t\,\wedge\,\zz^t\;+\;2\,t\,\nz^t\,\wedge\,(\,i^{N+1}\nz^t\,)
   \\
   \frac d{dt}\;\nn^t
   &=&
   -\;2\,t\,\nn^t\,\wedge\,\nn^t\;-\;2\,t\,\nz^t\,\wedge\,(\,i^{N+1}\nz^t\,)
  \end{eqnarray*}
 \end{Lemma}

 \proof
 Concerning the ordinary differential equation considered it may be
 interesting to know that the modified square $\nz\wedge(\,i^{N+1}\nz\,)$
 of an odd form $\nz$ is always concentrated in degrees divisible by $4$
 without constant term, because the number operator is additive on products:
 $$
  i^N(\;\nz\,\wedge\,(\,i^{N+1}\nz\,)\;)
  \;\;=\;\;
  (\,i^N\nz\,)\,\wedge\,(\,i^N\,i^{N+1}\nz\,)
  \;\;=\;\;
  (\,i^N\nz\,)\,\wedge\,(\,-\,i\,\nz\,)
  \;\;=\;\;
  \nz\,\wedge\,(\,i^{N+1}\nz\,)
 $$
 In consequence the standard theorems about the existence and uniqueness
 of solutions to ordinary differential equations guarantee the existence
 for all times of a unique solution staying in $\L^\odd T^*\oplus
 \L^{4\bullet}_\circ T^*\oplus\L^{4\bullet}_\circ T^*$ for every initial
 value in $\L^\odd T^*\oplus\L^{4\bullet}_\circ T^*\oplus\L^{4\bullet}_\circ
 T^*$. Of course for the ordinary differential equation considered a
 straightforward induction on degree shows that the homogeneous components
 of $\nz^t$, $\zz^t$ and $\nn^t$ of degree $r$ are polynomials in $t$ and
 the homogeneous components of degree less than $r$ of the initial values
 $\nz$, $\zz$ and $\nn$.

 On the other hand the power series $(\nz^t,\,\zz^t,\,\nn^t)$ in question
 certainly evaluates to the given intial value $(\nz^0,\,\zz^0,\,\nn^0)\,=\,
 (\nz,\,\zz,\,\nn)$ in $t=0$. In order to verify it is actually a solution
 to the ordinary differential equation considered we observe that the
 components $\zz^t$ and $\nn^t$ of the power series can be rewritten in
 terms of the component $\nz^t$ to read:
 \begin{eqnarray*}
  \zz^t
  &:=&
  \zz\,\wedge\,(\;1\;+\;t^2\,\zz\;)^{-1}
  \;+\;t^2\,\nz^t\,\wedge\,(\,i^{N+1}\nz\,)
   \,\wedge\,(\;1\;+\;t^2\,\zz\;)^{-1}
  \\
  \nn^t
  &:=&
  \nn\,\wedge\,(\;1\;+\;t^2\,\nn\;)^{-1}
  \;-\;t^2\,\nz^t\,\wedge\,(\,i^{N+1}\nz\,)
   \,\wedge\,(\;1\;+\;t^2\,\nn\;)^{-1}
 \end{eqnarray*}
 With $\nz^t$ being an odd form we know of course $\nz^t\,\wedge\,\nz^t\,=\,0$
 and thus we find immediately:
 $$
  \frac d{dt}\;\nz^t
  \;\;=\;\;
  \nz^t\,\wedge\,\Big(-2\,t\,\zz\,(\,1\,+\,t^2\,\zz\,)^{-1}
  \,-2\,t\,\nn\,(\,1\,+\,t^2\,\nn\,)^{-1}\Big)
  \;\;=\;\;
  -\,2\,t\,\nz^t\,\wedge\,(\,\zz^t\,+\,\nn^t\,)
 $$
 Verifying the differential equation for $\zz^t$ and $\nn^t$ is somewhat more
 involved. To begin with we rewrite the definition of the component $\zz^t$ of
 the power series $(\nz^t,\,\zz^t,\,\nn^t)$ a second time
 $$
  \zz^t
  \;\;:=\;\;
  \Big(\;\zz\;+\;t^2\,\nz^t\,\wedge\,(\,i^{N+1}\nz\,)\;\Big)
  \,\wedge\,(\;1\;+\;t^2\,\zz\;)^{-1}
 $$
 in order to use the product rule and the formula just derived for
 $\frac d{dt}\nz^t$ in the calculation:
 \begin{eqnarray*}
  \frac d{dt}\;\zz^t
  &=&
  \Big(\,2\,t\,\nz^t\,+\,t^2\,\nz^t\,\wedge\,
  \Big(\,-\,2\,t\,\zz\,\wedge\,(\,1\,+\,t^2\,\zz\,)^{-1}
  \,-\,2\,t\,\nn\,\wedge\,(\,1\,+\,t^2\,\nn\,)^{-1}\,\Big)\,\Big)
  \\
  &&
  \qquad\qquad\wedge\,
  (\,i^{N+1}\nz\,)\,\wedge\,(\,1\,+\,t^2\,\zz\,)^{-1}
  \;+\;\zz^t\,\wedge\,
  (\,-\,2\,t\,\zz\,\wedge\,(\,1\,+\,t^2\,\zz\,)^{-1}\,)
  \\[2pt]
  &=&
  -\,2\,t\,\zz^t\,\wedge\,\zz^t\;+\;2\,t\,\nz^t\,\wedge\,(\,i^{N+1}\nz\,)
  \,\wedge\,(\,1\,+\,t^2\,\zz\,)^{-1}
  \\
  &&
  \qquad\qquad\wedge\,\Big(\;1\,-\,t^2\,\zz\,\wedge\,(\,1\,+\,t^2\,\zz\,)^{-1}
  \,-\,t^2\,\nn\,\wedge\,(\,1\,+\,t^2\,\nn\,)^{-1}\,+\,t^2\,\zz^t\;\Big)
  \\
  &=&
  -\,2\,t\,\zz^t\,\wedge\,\zz^t\;+\;2\,t\,\nz^t\,\wedge\,(\,i^{N+1}\nz\,)
  \,\wedge\,(\,1\,+\,t^2\,\zz\,)^{-1}\,\wedge\,(\,1\,+\,t^2\,\nn\,)^{-1}
 \end{eqnarray*}
 In the last line we used the tautology $1\,-\,t^2\,\nn\,\wedge(\,1\,+\,t^2
 \,\nn\,)^{-1}\,=\,(\,1\,+\,t^2\,\nn\,)^{-1}$ and the fact that wedging with
 the odd form $\nz^t$ kills the difference $\zz^t\,-\,\zz\,\wedge\,(\,1\,+\,
 t^2\,\zz\,)^{-1}$. On the other hand the operator $i^{N+1}$ commutes by its
 very definition with taking the product
 $$
  i^{N+1}(\;\nz^t\;)
  \;\;=\;\;
  (\,i^{N+1}\nz\,)
  \,\wedge\,(\,1\,+\,t^2\,\zz\,)^{-1}
  \,\wedge\,(\,1\,+\,t^2\,\nn\,)^{-1}
 $$
 with a form concentrated in degrees divisible by $4$. Mutatis mutandis
 the same arguments imply the validity of the ordinary differential equation
 claimed for the component $\nn^t$.
 \qed

 \pfill
 Applying the solution formula of Lemma \ref{etf} in every point $p\,\in\,M$
 we conclude that the time dependent geometric differential forms $\nz^t$,
 $\zz^t$ and $\nn^t$ can be expressed completely in terms of the basic
 geometric differential forms $\nz$, $\zz$ and $\nn$ of a hypersurface
 $M\,\subset\,M^\out$, say
 \begin{eqnarray}
  \lefteqn{\sum_{k\geq 0}g(\,d^\nabla\ff,\,(R^t)^k\,\ff\,)}
  &&
  \label{nzf}
  \\
  &=&
  \Big(\,\sum_{k\geq 0}g(\,d^\nabla\ff,\,R^k\,\ff\,)\,\Big)
  \wedge\Big(\,1+t^2\!
  \sum_{{\scriptstyle k\geq 0}\atop{\scriptstyle k\;\odd}}
  g(\,\ff,\,R^k\,\ff\,)\,\Big)^{-1}
  \!\!\!\wedge\Big(\,1+t^2\!
  \sum_{{\scriptstyle k\geq 0}\atop{\scriptstyle k\;\mathrm{even}}}
  g(\,d^\nabla\ff,\,R^k\,d^\nabla\ff\,)\,\Big)^{-1}
  \nonumber
 \end{eqnarray}
 for the time dependent form $\nz^t$. What links this formula (\ref{nzf})
 to the logarithmic transgression form $\delta(\log\,F)(\left.TM^\out
 \right|_M,\nabla^t,\ff\wedge N)$ of a hypersurface $M\,\subset\,M^\out$
 associated to a multiplicative sequence of Pontryagin forms is that the
 traces occuring in Definition \ref{thyp}
 \begin{eqnarray}
  \tr\left(\frac{\ff\wedge N}{2\pi i}\,\Big(\frac{R^t}{2\pi i}
  \Big)^{2k-1}\right)
  &=&
  \frac1{(\,2\,\pi\,i\,)^{2k}}\;
  \tr\Big(\,\left(\,\ff^\#\otimes N\,-\,N^\#\otimes\ff\,\right)
  \,(\,R^t\,)^{2k-1}\,\Big)
  \nonumber
  \\
  &=&
  \frac{(\,-1\,)^k}{(\,4\pi^2\,)^k}\;\Big(\;
  g(\;\ff,\,(\,R^t\,)^{2k-1}\,N\;)\;-\;g(\;N,\,(\,R^t\,)^{2k-1}\,\ff\;)\;\Big)
  \nonumber
  \\
  &=&
  2\,t\,\frac{(\,-\,1\,)^{k-1}}{(\,4\,\pi^2\,)^k}\;
  g(\;d^\nabla\ff,\,(\,R^t\,)^{2k-2}\,\ff\;)
  \label{inis}
 \end{eqnarray}
 of the logarithmic transgression form are essentially the homogeneous
 components $\nz^t_{4k-1}$ of degree $4k-1$ of the time dependent geometric
 differential form $\nz^t$!
 
 \begin{Corollary}[Logarithmic Transgression Form of Permeable Hypersurfaces]
 \hfill\label{logt}\break
  Consider the linear interpolation $t\longmapsto\nabla^t:=\nabla+t\,
  \ff\wedge N$ between the two natural metric connections $\nabla^0:=\nabla$
  and $\nabla^1:=\left.\nabla^\out\right|_M$ on the restriction $\left.TM^\out
  \right|_M$ of the outer tangent bundle to a permeable hypersurface
  $M\subset M^\out$. The logarithmic transgression form for this interpolating
  curve and a multiplicative sequence parametrized by $F(z)=1+O(z^2)$ reads
  $$
   \delta\,(\log\,F)(\;\left.TM^\out\right|_M,\,\nabla^t,\,\ff\wedge N\;)
   \;\;=\;\;
   \frac{d}{dt}\,\left(\;\frac{t^2\,\nz}{4\pi^2}\,
   LF\Big(\,\sqrt{\frac{t^2\,d\nz}{4\pi^2}}\,\Big)\;\right)
  $$
  in terms of the $3$--form $\nz:=g(d^\nabla\ff,\ff)$ and the power series
  $LF(z):=\frac{\log\,F(z)}{z^2}$, in consequence:
  $$
   (\,\mathrm{trans}\,F\,)
   (\,\left.TM^\out\right|_M,\nabla,\nabla+\ff\wedge N\,)
   \;\;=\;\;
   \int_0^1
   F\Big(\,\sqrt{\frac{t^2\,d\nz}{4\pi^2}}\,\Big)
   \;\frac{d}{dt}\left(\;\frac{t^2\,\nz}{4\pi^2}\,
   LF\Big(\,\sqrt{\frac{t^2\,d\nz}{4\pi^2}}\,\Big)\;\right)\,dt
  $$
  With $F$ and $LF$ being even power series the ill--defined square root
  $\sqrt{d\nz}$ never materializes.
 \end{Corollary}

 \proof
 Recall that a permeable hypersurface is characterized by the vanishing
 of the $3$--form $R\,\ff=0$. In consequence the basic geometric form
 $\zz\,=\,0$ of a permeable hypersurface $M\,\subset\,M^\out$ vanishes
 identically, while $\nz\,=\,\nz_3\,=\,g(\,d^\nabla\ff,\ff\,)$ is a pure
 $3$--form with exterior derivative $d\nz\,=\,\nn\,=\,\nn_4$. The expansion
 (\ref{nzf}) of the time dependent differential form $\nz^t$ thus simplifies
 to $\nz^t\,=\,\nz\,\wedge\,(\,1\,+\,t^2\,d\nz\,)^{-1}$ so that the formula
 (\ref{inis}) for the traces occurring in the definition of the logarithmic
 transgression form turns into the formula
 \begin{eqnarray*}
  \delta\,(\,\log\,F\,)(\;\left.TM^\out\right|_M,\,\nabla^t,\,\ff\wedge N\;)
  &=&
  \sum_{k>0}k\,f_k\,\tr\left(\frac{\ff\wedge N}{2\pi i}\,
  \Big(\frac{R^t}{2\pi i}\Big)^{2k-1}\right)
  \\
  &=&
  \frac{d}{dt}\,\left(\;\sum_{k>0}f_k\,\frac{t^2\,\nz}{4\pi^2}\,\wedge\,
  \Big(\,\frac{t^2\,d\nz}{4\pi^2}\,\Big)^{k-1}\;\right)
 \end{eqnarray*}
 where the $f_k,\,k>0,$ are the coefficients of the power series
 $(\log\,F)(z)\,=:\,\sum_{k>0}f_kz^{2k}$ as before. Applying the
 exterior derivative to the logarithmic transgression form we get
 $$
  d\,\delta\,(\log\,F)
  (\left.TM^\out\right|_M,\nabla^\tau,\ff\wedge N)
  \;\;=\;\;
  \frac{d}{d\tau}
  \left(\sum_{k>0}f_k\Big(\frac{\tau^2\,d\nz}{4\pi^2}\Big)^k\right)
  \;\;=\;\;
  \frac{d}{d\tau}(\log\,F)\Big(\sqrt{\frac{\tau^2\,d\nz}{4\pi^2}}\Big)
 $$
 the integrand of the interior integral in the definition (\ref{deft}) of
 the transgression form:
 \begin{eqnarray*}
  \lefteqn{(\,\mathrm{trans}\,F\,)
  (\,\left.TM^\out\right|_M,\nabla,\nabla+\ff\wedge N\,)}
  \qquad&&
  \\
  &=&
  \int_0^1\exp\left(\;\int_0^t
  \frac{d}{d\tau}(\log\,F)\Big(\sqrt{\frac{\tau^2\,d\nz}{4\pi^2}}\Big)
  \,d\tau\;\right)\;
  \frac{d}{dt}\,\left(\;\frac{t^2\,\nz}{4\pi^2}\,
  LF\Big(\,\sqrt{\frac{t^2\,d\nz}{4\pi^2}}\,\Big)\;\right)\;dt
  \hbox to0pt{\quad$\;\;\qed$\hss}
 \end{eqnarray*}
\section{Special Values of Multiplicative Sequences}\label{bmsq}
 In this section we will use the transgression forms of hypersurfaces
 discussed in Section \ref{htra} in order to calculate some special values
 of multiplicative sequences of Pontryagin forms introduced in Section
 \ref{trag}. After a brief discussion of the Berger metrics $g^\rho$ of
 parameter $\rho\,>\,-1$ on the unit sphere $S^{2n-1}\,\subset\,\p$ of
 a hermitean vector space $\p$ we calculate the integral of the differential
 form $F(TB^{2n},\nabla^{\col,\rho})$ over the unit ball $B^{2n}$ for the
 Levi--Civita connection $\nabla^{\col,\rho}$ of a suitable collar metric
 $g^{\col,\rho}$, which induces to a Berger metric $g^\rho$ on the boundary
 $S^{2n-1}\,=\,\partial B^{2n}$. The collar metric $g^{\col,\rho}$ used in
 our calculations is actually the metric of the symmetric space $\C P^n$ and
 so we will rely heavily on the classical description of the Riemannian
 geometry of a symmetric space in terms of its associated triple product.

 \pfill
 Let us begin our calculations by a discussion of the natural CR--structure
 on odd dimensional spheres $S^{2n-1},\,n\geq\,1$. Consider a hermitean vector
 space $\p$ of complex dimension $n$. The real part $g$ of the hermitean
 form is a scalar product on $\p$ making $\p$ a euclidian vector space of
 dimension $2n$ endowed with an orthogonal complex structure $I\,\in\,\End\,
 \p$ satisfying $I^2\,=\,-1$ and $g(IX,IY)\,=\,g(X,Y)$ for all $X,\,Y\,\in\,
 \p$. The identification $T_P\p\,\cong\,\p$ of the tangent spaces in points
 $P\,\in\,\p$ with $\p$ itself makes $g$ a (flat) K\"ahler metric with
 integrable complex structure $I$ on $\p$ with associated K\"ahler form
 $\omega(X,Y)\,:=\,g(IX,Y)$ and Riemannian volume form $\frac1{n!}\,\omega^n$
 for the induced orientation. The K\"ahler metric $g$ restricts to the
 standard round metric on the unit sphere $S^{2n-1}\,\subset\,\p$,
 moreover the restrictions of the vector field and $1$--form
 $$
  C_P\;\;:=\;\;IP
  \qquad\qquad
  \gamma_P(\,X\,)\;\;:=\;\;C^\#_P(\,X\,)\;\;=\;\;g(\;IP,\,X\;)
 $$
 on $\p$ to $S^{2n-1}$ are called the Reeb vector field and the contact form
 respectively. Evidently the trivial connection $\nabla^{\mathrm{triv}}$ on
 $\p$ is torsion free and thus the Levi--Civita connection for the K\"ahler
 metric $g$, in particular we find $d\gamma\,=\,2\,\omega$. The Euler vector
 field $E_P\,:=\,P$ on $\p$ restricts to the outward pointing unit normal
 field on $S^{2n-1}$ and satisfies $E\ins\omega\,=\,\gamma$, hence the
 Riemannian volume form the standard round metric $g$ on $S^{2n-1}\,
 \subset\,\p$ reads:
 $$
  \vol_g
  \;\;:=\;\;
  (\,E\,\ins\,\frac1{n!}\,\omega^n\,)\Big|_{S^{2n-1}}
  \;\;=\;\;
  \frac1{2^{n-1}\,(n-1)!}\,\gamma\,\wedge\,(\,d\gamma\,)^{n-1}
 $$
 In consequence the restriction of $\gamma$ to $S^{2n-1}$ is really a contact
 form on $S^{2n-1}$, moreover
 \begin{equation}\label{voli}
  \<\,\frac{\gamma}{2\pi}\,\wedge\,(\,\frac{d\gamma}{2\pi}\,)^{n-1},
  \,[\,S^{2n-1}\,]\,>
  \;\;=\;\;
  \frac{\Gamma(n)}{2\,\pi^n}\int_{S^{2n-1}}
  \frac1{2^{n-1}(n-1)!}\,\gamma\,\wedge\,(\,d\gamma\,)^{n-1}
  \;\;=\;\;
  1
 \end{equation}
 due to the standard volume $\frac{2\,\pi^n}{\Gamma(n)}$ of the unit
 sphere in dimension $2n-1$. In general the Reeb vector field on a
 contact manifold with contact form $\gamma$ is the unique vector field
 $C$ satisfying $\gamma(C)\,=\,1$ and $C\ins d\gamma\,=\,0$, properties
 evidently satisfied by our Reeb field $C$.

 Digressing for a second we recall that every Killing vector field $K$ on a
 Riemannian manifold $M$ with Levi--Civita connection $\nabla$ and curvature
 tensor $R$ determines a skew symmetric endomorphism field $\mathfrak{K}\,
 \in\,\Gamma(\L^2TM)$ such that $(\,K,\,\mathfrak{K}\,)$ satisfies the
 extended Killing equation:
 $$
  \nabla_XK\;\;=\;\;\mathfrak{K}\,X
  \qquad\qquad
  \nabla_X\mathfrak{K}\;\;=\;\;R_{X,\,K}
 $$
 In fact a trilinear form skew symmetric in two arguments is uniquely
 determined by its skew symmetrization in some other two arguments.
 Both trilinear forms $g(\,R_{X,K}Y,\,Z\,)$ and $g(\,(\nabla_X\mathfrak{K})Y,
 \,Z\,)\,=\,g(\,\nabla^2_{X,Y}K,\,Z\,)$ are skew symmetric in $Y$ and $Z$
 however and share the same skew symmetrization $g(\,R_{X,Y}K,\,Z\,)$ in
 $X$ and $Y$ by the first Bianchi identity.

 The vector field $C$ provides a nice example of the extended Killing
 equation. Evidently $C$ is a Killing field for the flat K\"ahler
 metric $g$ on $\p$ with parallel skew symmetric covariant derivative
 $\nabla^{\mathrm{triv}}C\,=\,I$ consistent with $R^{\mathrm{triv}}\,=\,0$.
 After restriction to the unit sphere $S^{2n-1}$ the same argument becomes
 more interesting, the Reeb field $C$ is a Killing field on $S^{2n-1}$
 for the standard round metric $g$ with curvature $R_{X,Y}\,=\,-\,
 (X\wedge Y)$, hence:
 \begin{equation}\label{derc}
  \nabla_XC\;\;=\;\;\K\,X
  \qquad\qquad
  \nabla_X\K\;\;=\;\;C\,\wedge\,X
 \end{equation}
 where the skew symmetric endomorphism field $\K$ evaluates in $P\in S^{2n-1}$
 to the restriction
 $$
  \K_P\,X
  \;\;:=\;\;
  \pr_{\{\,P\,\}^\perp}(\;\nabla^{\mathrm{triv}}_XC\;)
  \;\;=\;\;
  \left\{\begin{array}{cl}
   IX & \textrm{for\ }X\,\in\,\{\,C_P\,\}^\perp
   \\[2pt]
   0  & \textrm{for\ }X\,\in\,\R C_P
  \end{array}\right.
 $$
 of the complex structure $I$ to $T_PS^{2n-1}\,=\,\{\,P\,\}^\perp$ defining
 the CR--structure on $S^{2n-1}$. Concluding this brief introduction to the
 CR--geometry of odd dimensional spheres we define the standard Berger metric
 $g^\rho$ with parameter $\rho\,>\,-1$ by scaling the radius of Hopf circles,
 the intersections of $S^{2n-1}\,\subset\,\p$ with complex lines in $\p$, by
 the factor $\sqrt{1+\rho}$
 \begin{equation}\label{berger}
  g^\rho\;\;:=\;\;g\;+\;\rho\,\gamma\,\otimes\,\gamma
 \end{equation}
 every Riemannian metric proportional to $g^\rho$ will be called a Berger
 metric of parameter $\rho$. Unless otherwise stated however the musical
 isomorphisms $\#$ and $\b$ will always refer to the standard metric $g$
 on $S^{2n-1}$ regardless of which Berger metric we are studying. For
 the Levi--Civita connection of $g$ we find say $\nabla_X\gamma\,=\,
 (\nabla_XC)^\#\,=\,g(\,\K X,\cdot\,)$ using (\ref{derc}) and obtain
 \begin{eqnarray*}
  (\,\nabla_Xg^\rho\,)(\;Y,\,Z\;)
  &=&
  \rho\;g(\;\K X,\,Y\;)\,\gamma(\,Z\,)\;+\;
  \rho\;\gamma(\,Y\,)\,g(\;\K X,\,Z\;)
  \\
  &=&
  g^\rho(\;Y,\,\rho\;(\,\gamma\cdot\K\,)_XZ\;)\;+\;
  g^\rho(\;\rho\;(\,\gamma\cdot\K\,)_XY,\,Z\;)
 \end{eqnarray*}
 where the endomorphism--valued $1$--form $\gamma\cdot\K$ on $S^{2n-1}$ is
 defined via:
 $$
  (\,\gamma\,\cdot\,\K\,)_X
  \;\;:=\;\;
  \gamma\otimes\K X\;+\;\gamma(\,X\,)\,\K
 $$
 Since $(\gamma\cdot\K)_XY$ is symmetric in $X$ and $Y$ we conclude that
 $\nabla^\rho\,:=\,\nabla\,+\,\rho\,(\gamma\cdot\K)$ defines the Levi--Civita
 connection for every metric proportional to the standard Berger metric
 $g^\rho$.

 \pfill
 Geometrically the Berger metrics $g^\rho,\,\rho>-1,$ arise as the metrics
 on $S^{2n-1}$ induced on the distance spheres in the complex projective
 spaces. With the complex projective spaces being in particular symmetric
 spaces it seems prudent to make use of the classical canon of formulas
 \cite{helg} relating the geometry of a symmetric space to its associated
 triple product. Consider therefore the involutive automorphism of the group
 $\SU(n+1)$ of special unitary transformations of $\C^{n+1}$ given by
 conjugation $\theta:\,\SU(n+1)\longrightarrow\SU(n+1),\,A\longmapsto
 SAS,$ with the reflection $S$ along the first standard basis vector. The
 automorphism $\theta$ fixes the subgroup $\mathbf{S}(\U(1)\times\U(n))
 \,\subset\,\SU(n+1)$ of special unitary transformations preserving the
 orthogonal decomposition $\C^{n+1}\,=\,\C\oplus\C^n$, hence the quotient
 $\SU(n+1)/\mathbf{S}(\U(1)\times\U(n))$ can be interpreted as the set
 $\C P^n$ of all such decompositions of $\C^{n+1}$ into a line and its
 orthogonal complement. On the Lie algebra of trace--free skew hermitean
 matrices $\su(n+1)$ of $\SU(n+1)$ the eigenspace decomposition $\su(n+1)
 \,=\,\k\oplus\p$ under the involution $\theta$ reads:
 $$
  \k\;\;:=\;\;\{\;\pmatrix{a&0\cr 0&A}\;|\;A\,=\,-A^H,\;a\,+\,\tr\,A\,=\,0\;\}
  \qquad
  \p\;\;:=\;\;\{\;\pmatrix{0&-X^H\cr X&0}\;|\;X\,\in\,\C^n\;\}
 $$
 The characterizing property of the Fubini--Study metric $g^\FS$ on $\C P^n$
 is that the corresponding scalar product on $\p$ is the standard scalar
 product $g(X,Y)\,:=\,\mathrm{Re}\;X^HY$ arising from the isomorphism
 $\p\,\cong\,\C^n$ indicated above, evidently $g$ is the restriction of
 the $\SU(n+1)$--invariant trace form $B(X,Y)\,=\,-\frac12\,\mathrm{Re}
 \;\tr\,XY$ on $\su(n+1)$ to $\p$. The triple product of the complex projective
 space $\C P^n$ is then defined on its isotropy representation $\p\,\cong\,
 \C^n$ via:
 \begin{eqnarray*}
  \lefteqn{[\;{\textstyle\frac12}\,X^2,\,Y\;]}\quad
  &&
  \\
  &:=&
  \Big[\;\pmatrix{0&-X^H\cr X&0},\;\Big[\;\pmatrix{0&-X^H\cr X&0},
  \;\pmatrix{0&-Y^H\cr Y&0}\;\Big]\;\Big]
  \\[3pt]
  &=&
  \pmatrix{0&+\,X^HXY^H\,-\,2\,X^HYX^H\,+\,Y^HXX^H
  \cr-\,XX^HY\,+\,2\,XY^HX\,-\,YX^HX&0}
 \end{eqnarray*}
 Alternatively this triple product can be written using only the scalar
 product $g$
 \begin{equation}\label{tpc}
  [\;{\textstyle\frac12}\,X^2,\,Y\;]
  \;\;=\;\;
  g(\,X,\,Y\,)\,X\;-\;3\,g(\,IX,\,IY\,)\,IX\;-\;g(\,X,\,X\,)\,Y
 \end{equation}
 and the canonical complex structure on $\p$ considered as an endomorphism
 $I$ on the real vector space underlying $\p\,=\,\C^n$. As said before the
 triple product associated to a symmetric space describes its geometry more
 or less completely by means of the following theorem \cite{helg}:

 \begin{Theorem}[Symmetric Metric in Conical Exponential Coordinates]
 \hfill\label{corm}
  Consider the isotropy representation $\p$ of a Riemannian symmetric space
  $G/K$ as a euclidian vector space of dimension $m$ with scalar product $g$,
  unit sphere $S^{m-1}$ and a triple product $[\,,]$. In consistency with the
  Lemma of Gau\ss\ the Riemannian metric on $G/K$ pulls back under conical
  exponential coordinates $\R^+\times S^{m-1}\longrightarrow G/K,\,(r,P)
  \longmapsto\exp(rP)K,$ to the metric
  $$
   {\widetilde\exp}^*g^{G/K}\;\;=\;\;dr\,\otimes\,dr\;+\;g^{G/K}_r
  $$
  where $g^{G/K}_r$ for sufficiently small $r>0$ is the Riemannian metric
  on $S^{m-1}_r\,:=\,\{r\}\times S^{m-1}$ defined in $P\in S^{m-1}_r$ in
  terms of the endomorphism $(\ad\;rP)^2\,:=\,r^2\,[\frac12\,P^2,\cdot]$
  by the formula:
  $$
   (\,g^{G/K}_r\,)_P(\;X,\,Y\;)
   \;\;=\;\;r^2\,
   g(\;\frac{\sinh\;\ad\;rP}{\ad\;rP}\,X,\;\frac{\sinh\;\ad\;rP}{\ad\;rP}\,Y\;)
  $$
 \end{Theorem}

 \noindent
 For the specific triple product (\ref{tpc}) the endomorphisms $(\ad\,P)^2:
 \,\p\longrightarrow\p,\,X\longmapsto[\,\frac12P^2,\,X\,],$ are evidently
 symmetric for all unit vectors $P\,\in\,S^{2n-1}\,\subset\,\p$ with
 eigenvalues $0$, $-4$ and $-1$ on the pairwise orthogonal eigenspaces
 $\R\,P$, $\R\,IP$ and $\{P,IP\}^\perp$. Formally at least the root
 endomorphisms $\ad\,P$ have eigenvalues $0$, $2i$ and $i$, although
 a minor nuisance the only effect of this appearance of $i$ is to convert
 the hyperbolic sines of Theorem \ref{corm} via $\frac{\sinh\,i\lambda}
 {i\lambda}\,=\,\frac{\sin\,\lambda}\lambda$ into sines. In consequence
 the Fubini--Study metric pulls back under conical exponential coordinates
 $\widetilde\exp:\,\R^+\times S^{2n-1}\longrightarrow\C P^n,\,(r,P)\longmapsto
 \exp(rP)\,K,$ to a smooth metric on the open ball $B_{\frac\pi2}(0)\,\subset
 \,\p$ of radius equal to the injectivity radius $\frac\pi2$ of $\C P^n$
 \begin{equation}\label{fsm}
  \widetilde\exp^*\,g^\FS\;\;=\;\;
  dr\,\otimes\,dr\;+\;\sin^2\,r\,
  (\;g\;-\;\sin^2r\,(\,\gamma\,\otimes\,\gamma\,)\;)
 \end{equation}
 where we have used $(r\frac{\sin\,2r}{2r})^2\,=\,\sin^2r\,(\,1-\sin^2r\,)$.
 In particular the Fubini--Study metric $\widetilde\exp^*g^\FS$ restricts on
 the distance spheres $S^{2n-1}_r$ of radius $0<r<\frac\pi2$ to the Berger
 metric $-\rho\,g^\rho$ with parameter $\rho\,:=\,-\sin^2r\,\in\,]-1,\,0\,[$.
 In passing we remark that the considerations leading to Theorem \ref{corm}
 can be repeated with every parallel tensor on an arbitrary symmetric space
 $G/K$, mimicking in particular the calculations for the Riemannian metric
 $g^\FS$ we obtain the following formula for the pull back of the parallel
 K\"ahler form $\omega^\FS$ on $\C P^n$:
 \begin{equation}\label{fso}
  {\widetilde\exp}^*\,\omega^\FS\;\;=\;\;
  \sin\,r\;\cos\,r\;dr\,\wedge\,\gamma\;+\;\frac{\sin^2r}2\;d\gamma
 \end{equation}
 In this and all subsequent formulas $\gamma$ and $d\gamma$ denote the pull
 back to $\R^+\times S^{2n-1}$ of the differential forms of the same name on
 the unit sphere $S^{2n-1}\subset\p$, in particular $\frac\partial{\partial r}
 \ins d\gamma\,=\,0$.

 \pfill
 In essence formula (\ref{fsm}) provides us with a very specific Riemannian
 metric on the closed ball $\overline B_r(0)\,\subset\,\p$ of radius $0\,<\,
 r\,<\,\frac\pi2$, which induces a Berger metric $-\rho\,g^\rho$ with parameter
 $\rho\,=\,-\sin^2r\,\in\,]-1,\,0\,[$ on boundary sphere $S^{2n-1}_r$. All we
 still need to do in order to proceed is to calculate the second fundamental
 form of the boundary sphere $S^{2n-1}_r$. For this purpose we observe that
 the second fundamental form of an arbitrary hypersurface can be calculated
 as the Lie derivative $\ff\,=\,-\frac12\,\mathfrak{Lie}_Ng$ of the metric
 along the normal vector field $N$ used to define $\ff$ in the first place.
 For the boundary spheres $S^{2n-1}_r\,\subset\,\R^+\times S^{2n-1}$ the normal
 vector field of choice is the outward pointing normal field $\frac\partial
 {\partial r}$ so that we obtain
 \begin{eqnarray}
  \ff
  &=&
  -\,{\textstyle\frac12}\,\mathfrak{Lie}_{\frac\partial{\partial r}}
  \Big(\;dr\,\otimes\,dr\;+\;
  \sin^2r\,(\;g\;-\;\sin^2r\,(\,\gamma\,\otimes\,\gamma\,)\;)\;\Big)
  \label{ff}
  \\
  &=&
  -\,\sin\,r\,\cos\,r\,(\;g\;-\;2\,\sin^2r\,(\,\gamma\,\otimes\,\gamma\,)\;)
  \nonumber
 \end{eqnarray}
 for the second fundamental $\ff$ considered as a symmetric bilinear form
 on the tangent space of $S^{2n-1}_r$. Rewriting this result into a vector
 valued $1$--form on $S^{2n-1}_r$ we get in turn
 \begin{equation}\label{sff}
  \ff
  \;\;=\;\;
  -\,\frac{\cos\,r}{\sin\,r}\,\Big(\;(\id\,-\,\gamma\,C)
  \;+\;\frac{1-2\sin^2r}{1-\sin^2r}\,(\gamma\,C)\;\Big)
  \;\;=\;\;
  -\,\frac1{\tan\,r}\,\Big(\;\id\;+\;\frac\rho{\rho\,+\,1}\,\gamma\,C\;\Big)
 \end{equation}
 using $\rho\,:=\,-\sin^2r$ as before. Being torsion--free the Levi--Civita
 connection $\nabla^\rho$ for the Berger metric $-\rho g^\rho$ on $S^{2n-1}_r$
 kills the soldering form $d^{\nabla^\rho}\id\,=\,0$ while $\nabla^\rho C\,=\,
 (\rho+1)\K$, so
 $$
  d^{\nabla^\rho}\ff
  \;\;=\;\;
  \tan\,r\;\Big(\;d\gamma\,C\;-\;(\,\rho\,+\,1\,)\,\gamma\,\K\;\Big)
 $$
 where $\K$ has to be interpreted as a vector valued $1$--form of course.
 Recalling the definition $\nz\,:=\,-\rho g^\rho(d^{\nabla^\rho}\ff,\ff)$ of
 of the basic differential $3$--form $\nz$ of the distance sphere $S^{2n-1}_r$
 we find
 \begin{equation}\label{xif}
  \nz
  \;\;=\;\;
  \rho\,g^\rho\Big(\;d\gamma\,C\;-\;(\,\rho\,+\,1\,)\,\gamma\,\K,\,
  \id\,+\,\frac\rho{\rho\,+\,1}\,\gamma\,C\;\Big)
  \;\;=\;\;
  \rho^2\,\gamma\,\wedge\,d\gamma
 \end{equation}
 omitting the auxiliary calculations:
 $$
  g^\rho(\,C,\id\,) \;\,=\;\,(\,\rho\,+\,1\,)\,\gamma
  \quad\;\;
  g^\rho(\,C,C\,)\; \,=\;\,\rho+1
  \quad\;\;
  g^\rho(\,\K,\id\,)\;\,=\;\,d\gamma
  \quad\;\;
  g^\rho(\,\K,C\,)  \;\,=\;\,0
 $$
 Strictly speaking $\nz$ is only the homogeneous part of degree three of the
 form $\nz$ defined in Section \ref{htra}, however all the higher degree parts
 vanish, because $S^{2n-1}_r\,\subset\,\R^+\times S^{2n-1}$ is a permeable
 hypersurface! According to Definition \ref{perm} a permeable hypersurface
 is characterized by a covariantly closed $d^{\nabla^\rho}(d^{\nabla^\rho}\ff)
 \,=\,0$ second fundamental form. Due to the extended Killing equation
 (\ref{derc}) for $\nabla_X\K\,=\,C\wedge X$ the covariant derivative
 of $\K$ with respect to $\nabla^\rho$
 \begin{eqnarray*}
  (\,\nabla^\rho_X\K\,)\,Y
  &=&
  (\,\nabla_X\K\,)\,Y\;+\;\rho\,
  [\;\gamma(X)\,\K\;+\;\gamma\,\otimes\,\K X,\;\K\;]\,Y
  \\[2pt]
  &=&\gamma(Y)\,X\;-\;g(\,X,\,Y\,)\,C\;-\;\rho\,\gamma(Y)\,\K^2X
 \end{eqnarray*}
 becomes $d^{\nabla^\rho}\K\,=\,-\gamma(\id-\rho\K^2)$ after skew
 symmetrization in $X,\,Y$ so that $d^{\nabla^\rho}\id\,=\,0$ implies:
 $$
  d^{\nabla^\rho}(\;d^{\nabla^\rho}\ff\;)
  \;\;=\;\;
  \tan\,r\;\gamma\;d^{\nabla^\rho}(\,\nabla^\rho C\,)
  \;\;=\;\;
  -\,\tan\,r\;(\,\rho\,+\,1\,)\;\gamma^2\;(\,\id\,-\,\rho\,\K^2\,)
  \;\;=\;\;
  0
 $$
 By Corollary \ref{logt} the logarithmic transgression form of the distance
 sphere $S^{2n-1}_r$ reads
 \begin{equation}\label{logs}
  \delta\,(\log\,F)(\;\left.T(\R^+\times S^{2n-1})\right|_{S^{2n-1}_r},
  \,(\nabla^\rho)^t,\,\ff\wedge\frac\partial{\partial r}\;)
  \;\;=\;\;
  \frac{d}{dt}\,\left(\;\frac{t\,\rho\,\gamma}{2\pi}\,
  \frac{t\rho\,d\gamma}{2\pi}\,
  LF\Big(\,\frac{t\,\rho\,d\gamma}{2\pi}\,\Big)\;\right)
 \end{equation}
 for every multiplicative sequence of Pontryagin forms parametrized by an
 even power series $F(z)\,=\,1+O(z^2)$ with logarithm $LF(z)\,=\,\frac
 {\log\,F(z)}{z^2}$. With the decisive equations (\ref{creq}), (\ref{fso})
 and (\ref{logs}) in place we can eventually embark on the proof of the main
 theorem of this article:

 \begin{Theorem}[Special Values of Multiplicative Sequences]
 \hfill\label{null}\break
  Consider a smooth metric $g^{\col,\rho}$ on the closed unit ball
  $B^{2n}\,\subset\,\p$, which is a product of the usual metric on
  $]-\varepsilon,\,0\,]$ and the Berger metric of parameter $\rho\in\,]-1,0\,[$
  in a collar neighborhood $]-\varepsilon,\,0\,]\times S^{2n-1}$ of the
  boundary $S^{2n-1}$ of $B^{2n}$. The multiplicative sequence of Pontryagin
  forms parametrized by an even formal power series $F(z)=1+O(z^2)$ with
  associated composition inverse $\phi(z)$ of $\frac z{F(z)}=z+O(z^3)$
  takes the following value on $B^{2n}$:
  $$
   \int_{B^{2n}}F(\;TB^{2n},\,\nabla^{\col,\rho}\;)
   \;\;=\;\;
   \rho^n\,\mathrm{res}_{z=0}\left[\;\frac{F(\,z\,)^n}{z^{n+1}}\,dz\;\right]
   \;\;=\;\;
   \rho^n\,\mathrm{res}_{z=0}\left[\;\frac{(\log\,\phi)'(z)}{z^n}\,dz\;\right]
  $$
 \end{Theorem}

 \proof
 The special Berger metric $-\rho\,g^\rho$ with parameter $\rho\,\in\,
 ]-1,\,0\,[$ is realized by the geodesic distance sphere $S^{2n-1}_r\,
 \cong\,S^{2n-1}$ in $\C P^n$ of radius $r\,\in\,]\,0,\,\frac\pi2[$ with
 $\rho\,=\,-\sin^2r$. The principal idea of the proof is to transgress
 from the Levi--Civita connection $\nabla^{\col,\rho}$ for the collar
 metric $g^{\col,\rho}$ to the Levi--Civita connection $\nabla^\FS$ of
 the Fubini--Study metric $g^\FS$ of $\C P^n$ restricted to the ball
 $B^{2n}_r\,\cong\,B^{2n}_r\,\subset\,\p$ of radius $r$. For the collar
 metric $g^{\col,\rho}$ the boundary $S^{2n-1}\,=\,\partial B^{2n}_r$
 is totally geodesic, hence the transgression restricts on the boundary
 to the standard hypersurface transgression problem from $\nabla^\rho$
 to $\nabla^\FS$. For the Fubini--Study metric $g^\FS$ on the other hand
 the boundary $S^{2n-1}_r$ is permeable with the known logarithmic
 transgression form (\ref{logs}).

 In a first step we employ the naturality of the multiplicative sequence
 associated to the even formal power series $F$ under scaling $\mathrm{sc}:
 \,B^{2n}_r\longrightarrow B^{2n}$ and the transgression formula (\ref{teq})
 to convert the integral $\<F(TB^{2n},\nabla^{\col,\rho}),[B^{2n}]>\,=\,
 \<F(TB^{2n}_r,\mathrm{sc}^*\nabla^{\col,\rho}),[B^{2n}_r]>$ into:
 \begin{eqnarray}
  \lefteqn{\<\;F(\;TB^{2n}_r,\,\nabla^\FS\;),\,[\,B^{2n}_r\,]\;>
   \;-\;\<\;F(\;TB^{2n},\,\nabla^{\col,\rho}\;),\,[\,B^{2n}\,]\;>}
  \qquad\quad&&
  \label{fin}
  \\
  &=&\<\;F(\;TS^{2n-1}_r,\,\nabla^\rho\;)\,\wedge\,(\mathrm{trans}\,F)
  (\;\left.TB^{2n}_r\right|_{S^{2n-1}_r},\,\nabla^\rho,\,\nabla^\rho\,
  +\,\ff\,\wedge\,\frac\partial{\partial r}\;),\,[\,S^{2n-1}_r\,]\;>
  \nonumber
 \end{eqnarray}
 Recall from the discussion of complex projective spaces in Section \ref{trag}
 that $F(T\C P^n,\nabla^\FS)$ is parallel so that its top term is a constant
 multiple of the volume form of $\C P^n$. Using the formula (\ref{fsm}) for
 the Fubini--Study metric $g^\FS$ in conical exponential coordinates we find
 that
 \begin{eqnarray*}
  \mathrm{Vol}(\;B^{2n}_r,\,g^\FS\;)
  &=&
  \int_0^r\,\mathrm{Vol}(\;S^{2n-1},\,
  \sin^2s\,(\,g\,-\,\sin^2s\,\gamma\otimes\gamma\,)\;)\,ds
  \\
  &=&
  \mathrm{Vol}(\;S^{2n-1},\,g\;)
  \,\int_0^r\,\sqrt{1\,-\,\sin^2s}\;\sin^{2n-1}s\;ds
  \;\;=\;\;
  \frac{\pi^n}{n!}\,\sin^{2n}r
 \end{eqnarray*}
 extended to the injectivity radius $r\,=\,\frac\pi2$ of $\C P^n$ becomes
 its volume $\frac{\pi^n}{n!}$ and conclude
 \begin{equation}\label{first}
  \<\;F(\;TB^{2n}_r,\,\nabla^\FS\;),\,[\,B^{2n}_r\,]\;>
  \;\;=\;\;
  \rho^n\,\mathrm{res}_{z=0}\Big[\;\frac{F(z)^{n+1}}{z^{n+1}}\,dz\;\Big]
 \end{equation}
 from equation (\ref{creq}), the additional sign in $\rho\,=\,-\sin^2r$
 turns out to be immaterial, because both sides of the equation vanish if
 $n$ is odd.
 
 Moreover the explicit formula (\ref{fso}) for the K\"ahler form
 $\omega^\FS$ in conical exponential coordinates implies $\left.\omega^\FS
 \right|_{S^{2n-1}_r}\,=\,-\frac12\,\rho\,d\gamma$. Using $F(T\C P^n,
 \nabla^\FS)\,=\,F(\frac{\omega}\pi)^{n+1}$ we can rewrite the basic
 transgression formula (\ref{this}) for the logarithmic transgression
 form (\ref{logs}) as
 \begin{eqnarray*}
  \left.F(\,TB^{2n}_r,\,\nabla^\FS\,)\right|_{S^{2n-1}_r}
  &=&
  F(\,\left.\frac{\omega^\FS}\pi\right|_{S^{2n-1}_r}\,)^{n+1}
  \;\;=\;\;
  F(\;\frac{\rho\,d\gamma}{2\pi}\;)^{n+1}
  \\
  &=&
  F(\,TS^{2n-1}_r,\,\nabla^\rho\,)\wedge\left.\exp
  \Big(d\int_0^t\frac d{d\tau}\Big(\,\frac{\tau\rho\gamma}{2\pi}\,
  \frac{\tau\rho\,d\gamma}{2\pi}\,LF(\frac{\tau\rho\,d\gamma}{2\pi})
  \Big)d\tau\Big)\right|_{t=1}
  \\
  &=&
  F(\,TS^{2n-1}_r,\,\nabla^\rho\,)\wedge\left.\exp\Big(\;
  (\,\log\,F\,)(\,\frac{t\rho\,d\gamma}{2\pi}\,)\;\Big)\right|_{t=1}
 \end{eqnarray*}
 recall that $F$ is even with $LF(z)\,:=\,\frac{\log\,F(z)}{z^2}$. With the
 differential form $F(\frac{\rho d\gamma}{2\pi})$ being invertible in the
 algebra of differential forms we conclude $F(TS^{2n-1}_r,\nabla^\rho)\,=\,
 F(\frac{\rho d\gamma}{2\pi})^n$. Reinserting this result into the preceding
 equation we obtain for the transgression form in equation (\ref{fin})
 \begin{eqnarray*}
  \lefteqn{F(\;TS^{2n-1}_r,\,\nabla^\rho\;)\,\wedge\,(\mathrm{trans}\,F\,)
  (\;\left.TB^{2n}\right|_{S^{2n-1}_r},\,\nabla^\rho,\,\nabla^\rho\,+\,
  \ff^\rho\wedge N\;)}\qquad\qquad
  &&
  \\
  &=&
  F(\,\frac{\rho\,d\gamma}{2\pi}\,)^n
  \,\int_0^1F(\;\frac{t\rho\,d\gamma}{2\pi}\;)
  \,\frac d{dt}\Big(\;\frac{t\rho\gamma}{2\pi}
  \,\frac{t\rho\,d\gamma}{2\pi}\;LF(\;\frac{t\rho\,d\gamma}{2\pi}\;)\;\Big)\,dt
 \end{eqnarray*}
 which integrates via equation (\ref{voli}) in the form $\<\frac\gamma{2\pi}
 \wedge(\frac{d\gamma}{2\pi})^{n-1},[\,S^{2n-1}\,]>\,=\,1$ to the value
 \begin{eqnarray*}
  \lefteqn{\mathrm{res}_{z=0}\left[\,F(\rho z)^n\;\left(\int_0^1F(t\rho z)\,
  \frac d{dt}\Big(\,t^2\rho^2z\,\frac{\log\,F(t\rho z)}{(t\rho z)^2}\,\Big)
  \,dt\right)\frac{dz}{z^n}\,\right]}
  \;\;&&
  \\[3pt]
  &=&
  \mathrm{res}_{z=0}\left[\,F(\rho z)^n
  \left(\int_0^1F(t\rho z)\,\frac{F'(t\rho z)}{F(t\rho z)}\,\rho\,dt\right)
   \frac{dz}{z^n}\,\right]
  \;\;=\;\;
  \rho^n\,
  \mathrm{res}_{z=0}\left[\,\frac{F(z)^{n+1}-F(z)^n}{z^{n+1}}\,dz\,\right]
 \end{eqnarray*}
 with the formal variable $z\,:=\frac{d\gamma}{2\pi}$. Combined with
 equations (\ref{fin}) and (\ref{first}) this result implies
 $$
  \<\;F(\,TB^{2n},\,\nabla^{\col,\rho}\,),\,[\,B^{2n}\,]\;>
  \;\;=\;\;
  \mathrm{res}_{z=0}\Big[\,\frac{F(z)^n}{z^{n+1}}\,\Big]
  \;\;=\;\;
  \mathrm{res}_{z=0}\Big[\,\frac{\phi(z)^n\,\phi'(z)}{z^n\,\phi(z)^{n+1}}
  \,dz\,\Big]
 $$
 where the odd power series $\phi$ is the composition inverse of
 $\frac z{F(z)}\,=\,z+O(z^3)$ as before.
 \qed

 \pfill
 In order to complete the proof of Theorem \ref{mult} stated in the
 introduction we still have to argue that Theorem \ref{null} is not
 only true for the values $\rho\in\,]-1,\,0\,[$ of the Berger parameter
 $\rho$ realized by the distance spheres in complex projective space.
 Evidently it is sufficient to show that the value $\<F(\,TB^{2n},\,
 \nabla^{\col,\rho}\,),[\,B_{2n}\,]>$ of the multiplicative sequence
 parametrized by $F$ on the closed unit ball $B^{2n}$ with a collar
 metric $g^{\col,\rho}$ for the Berger metric $g^\rho$ on $S^{2n-1}
 \,=\,\partial B^{2n}$ is an analytic function in $\rho$. For this
 purpose we consider the auxiliary metric
 $$
  g^{\out,\rho}
  \;\;:=\;\;
  dr\,\otimes\,dr\;+\;r^2\,(\;g\;+\;\rho\,r^2\,\gamma\,\otimes\,\gamma\;)
 $$
 on $\R^+\times S^{2n-1}$, which is smooth in $0$ and positive definite
 on a neighborhood of the unit ball $B^{2n}\,\subset\,\p$. Clearly this
 auxiliary Riemannian metric induces the Berger metric $g^\rho$ on the
 boundary unit sphere $S^{2n-1}\,\subset\,B^{2n}$, nevertheless the boundary
 is not totally geodesic. Using the Lie derivative with respect to the unit
 normal field $\frac\partial{\partial r}$ as in equation (\ref{ff}) we get
 $$
  \ff
  \;\;=\;\;
  -\,\frac12\,\left.\mathfrak{Lie}_{\frac\partial{\partial r}}
  \Big(\;dr\,\otimes\,dr\;+\;r^2\,g\;+\;r^4\,\rho\,
  \gamma\,\otimes\,\gamma\;\Big)\right|_{r=1}
  \;\;=\;\;
  -\,g\;-\;2\,\rho\,\gamma\,\otimes\,\gamma
 $$
 or $\ff\,=\,-(\,\id\,+\,\frac\rho{\rho+1}\,\gamma\,C\,)$ expressed as a
 vector valued $1$--form. It is straightforward, but somewhat more work to
 calculate the Levi--Civita connection $\nabla^{\out,\rho}$ for the metric
 $g^{\out,\rho}$. Leaving the details of this calculation to the reader we
 only note that it is sufficient to do the calculations for $\frac\partial
 {\partial r}$ and vector fields $X,\,Y\,\in\,\Gamma(\,TS^{2n-1}\,)$ on
 $S^{2n-1}$ extended constantly in $r$ direction to $\R^+\times S^{2n-1}$.
 Eventually we find the formulas
 $$
  \nabla^{\out,\rho}_{\frac\partial{\partial r}}\frac\partial{\partial r}
  \;\;=\;\;
  0
  \qquad\qquad
  \nabla^{\out,\rho}_X\frac\partial{\partial r}
  \;\;=\;\;
  \nabla^{\out,\rho}_{\frac\partial{\partial r}}X
  \;\;=\;
  \frac1r\,X\;+\;\frac{\rho\,r}{\rho\,r^2\,+\,1}\,\gamma(X)\,C
 $$
 and:
 $$
  \nabla^{\out,\rho}_XY
  \;\;=\;\;
  \nabla_XY\,-\,r\,g(X,Y)\,\frac\partial{\partial r}
  \,+\,\rho\,\Big(\,r^2\,\gamma(X)\,\K Y\,+\,r^2\,\gamma(Y)\,\K X
  \,-\,2\,r\,\gamma(X)\,\gamma(Y)\,\frac\partial{\partial r}\,\Big)
 $$
 For $\rho\,=\,0$ these formulas evidently reduce to the well--known
 formulas for the trivial connection $\nabla^{\mathrm{triv}}$ on $\p$
 expressed in polar coordinates. More important is that $\nabla^{\out,\rho}$
 is a rational function of $\rho$ with a pole on the sphere of radius
 $\frac1{\sqrt{-\rho}}$. Hence the curvature tensor and the basic geometric
 forms on the boundary hypersurface $S^{2n-1}\,=\,\partial B^{2n}$ are all
 rational functions of $\rho$. The integrand differential forms occurring
 on the right hand side of the reformulation
 $$
  \int_{B^{2n}}\!F(\,TB^{2n},\nabla^{\col,\rho}\,)
  \;\;=\;\;
  \int_{B^{2n}}\!F(\,TB^{2n},\nabla^{\out,\rho}\,)
  \;-\;
  \int_{S^{2n-1}}\!\mathrm{Trans}\,F
  (\,TB^{2n},\nabla^{\col,\rho},\nabla^{\out,\rho}\,)
 $$
 of the transgression formula (\ref{tokes}) are thus rational functions
 of $\rho$ as well so that the left hand side $\<\,F(\,TB^{2n},\,
 \nabla^{\col,\rho}\,),\,[\,B^{2n}\,]\,>$ is certainly an analytic
 function of $\rho$.

\begin{thebibliography}{99999}
%
 \bibitem[A]{a}
  \textsc{Atiyah, M.~:}\quad
  \textit{K--Theory,\ }
  \textrm{Cambridge University Press 12 (1966).}
%
 \bibitem[APS]{aps}
  \textsc{Atiyah, M.~, Patodi, V.~\& Singer, I.~:}\quad
  \textit{Spectral asymmetry and Riemannian geometry I,\ }
  \textrm{Mathematical Proceedings of the Cambridge Philosophical
    Society 77 (1975), 43---69.}
%
 \bibitem[B\"a]{b}
  \textsc{B\"ar, Ch.~:}\quad
  \textit{Metrics with Harmonic Spinors,\ }
  \textrm{Geometric and Functional Analysis 6 (1996), 899---942.}
%
 \bibitem[Be]{bes}
  \textsc{Bechtluft--Sachs, St.~:}\quad
  \textit{The Computation of $\eta$--Invariants on Manifolds
    with Free Circle Action,\ }
  \textrm{Journal Functional Analysis 174 (2000), 251---263.}
%
 \bibitem[BGV]{bgv}
  \textsc{Berline, N.~Getzler, E.~\& Vergne, M.~:}\quad
  \textit{Heat Kernels and Dirac Operators,\ }
  \textrm{Grundlehren der mathematischen Wissenschaften 298, Springer (1992).}
%
 \bibitem[Hab]{hab}
  \textsc{Habel, M.~:}\quad
  \textit{Die $\eta$--Invariante der Berger Sph\"aren,\ }
  \textrm{Diplom thesis at the University of Hamburg (2000).}
%
 \bibitem[Hel]{helg}
  \textsc{Helgason, M.~:}\quad
  \textit{Differential Geometry, Lie Groups and Symmetric Spaces,\ }
  \textrm{Pure and Applied Mathematics 80, Academic Press (1978).}
%
 \bibitem[Hir]{hir}
  \textsc{Hirzebruch, F.~:}\quad
  \textit{Topological Methods in Algebraic Geometry,\ }
  \textrm{Grundlehren der mathematischen Wissenschaften 131, Springer (1966).}
%
 \bibitem[Hit1]{ht1}
  \textsc{Hitchin, N.~:}\quad
  \textit{Harmonic Spinors,\ }
  \textrm{Advances in Mathematics 14 (1975), 1---55.}
%
 \bibitem[Hit2]{ht2}
  \textsc{Hitchin, N.~:}\quad
  \textit{Einstein Metrics and the $\eta$--Invariant,\ }
  \textrm{Bolletino della Unione Matematica Italiana B (7) 11 (1997), 95---105.}
%
 \bibitem[K]{k}
  \textsc{Koh, D.~:}\quad
  \textit{The $\eta$--invariant of the Dirac operator on the Berger spheres,\ }
  \textrm{Diplom thesis at the University of Hamburg (2004).}
%
 \bibitem[LM]{lm}
  \textsc{Lawson, H.~B.~\& Michelsohn, L.~:}\quad
  \textit{Spin Geometry,\ }
  \textrm{Princeton Mathematical Series 38, Princeton University Press (1990).}
%
\end{thebibliography}
\end{document}